\documentclass[final]{siamltex}

\def\d{\delta} 
 
\def\e{{\epsilon}}

\def\norm#1{\|#1\|} 
\def\wb{{\bar w}}

\usepackage{algorithm}
\usepackage{algorithmic}
\usepackage{amsmath}
\usepackage{amssymb}
\usepackage{caption}
\usepackage{comment}
\usepackage{graphicx}
\usepackage{float}
\usepackage[latin1]{inputenc}
\usepackage[pdfborder={0 0 0}, colorlinks=true, linkcolor=blue, citecolor=red]{hyperref}\hypersetup{pdfborder=0 0 0}
\usepackage{mathrsfs}
\usepackage{multirow}
\usepackage{showkeys,cite}
\usepackage{subfigure}
\usepackage{wrapfig}
\usepackage{xcolor}
\usepackage[textsize=footnotesize]{todonotes}
\usepackage{ulem}
\graphicspath{{./plots/}}

\begin{document}

\title{
       eHDG: An Exponentially Convergent Iterative Solver for HDG Discretizations of Hyperbolic Partial
       Differential Equations. \thanks{This research was
    partially supported by DOE grants DE-SC0010518 and
    DE-SC0011118. We are grateful for the supports.}}


\author{Sriramkrishnan Muralikrishnan\footnotemark[1] \and Minh-Binh Tran\footnotemark[2] \and Tan Bui-Thanh\footnotemark[3]}

\footnotetext[1]{Department of Aerospace Engineering and Engineering Mechanics, The University of Texas at Austin, Austin, TX 78712,
    USA.}
\footnotetext[2]{Department of Mathematics, University of Wisconsin, Madison, WI 53706, USA.}

\footnotetext[3]{Department of Aerospace Engineering and Engineering Mechanics, and the Institute for Computational Engineering and Sciences, The University of Texas at Austin, Austin, TX 78712,
    USA.}

\bibliographystyle{siam}
\newcommand{\TODO}[1]{ \fbox{\parbox{3in}{\bf TODO: #1}}}

\newcommand{\grbf}[1] {\mbox{\boldmath${#1}$\unboldmath}}
\newcommand{\gbf}[1] {\mathbf{#1}}

\newcommand{\beq} {\begin{equation}}
\newcommand{\eeq} {\end{equation}}
\newcommand{\bdm} {\begin{displaymath}}
\newcommand{\edm} {\end{displaymath}}
\newcommand{\bit}{\begin{itemize}}
\newcommand{\eit}{\end{itemize}}
\newcommand{\bde}{\begin{description}}
\newcommand{\ede}{\end{description}}
\newcommand{\bce}{\begin{center}}
\newcommand{\ece}{\end{center}}
\newcommand{\ben} {\begin{enumerate}}
\newcommand{\een} {\end{enumerate}}
\newcommand{\bea} {\begin{eqnarray}}
\newcommand{\eea} {\end{eqnarray}}
\newcommand{\barr} {\begin{array}}
\newcommand{\earr} {\end{array}}
\newcommand{\bean} {\begin{eqnarray*}}
\newcommand{\eean} {\end{eqnarray*}}
\newcommand{\edoc} {\end{document}}
\newcommand{\On}{\hspace{0.2cm} \mbox{on} \hspace{0.2cm}}
\newcommand{\In}{\hspace{0.2cm} \mbox{in} \hspace{0.2cm}}
\newcommand{\Minimize}{\hspace{0.2cm} \mbox{minimize} \hspace{0.2cm}}
\newcommand{\Maximize}{\hspace{0.2cm} \mbox{maximize} \hspace{0.2cm}}
\newcommand{\SubjectTo}{\hspace{0.2cm} \mbox{subject} \hspace{0.15cm}
            \mbox{to} \hspace{0.2cm}}

\def\etal{{\it et al.~}}

\newcommand{\point}[1] {\ensuremath{\boldsymbol{#1}}}
\newcommand{\vvec}[1] {\ensuremath{\boldsymbol{#1}}}
\renewcommand{\vec}[1] {\ensuremath{\boldsymbol{#1}}}
\newcommand{\cvec}[1] {\ensuremath{\boldsymbol{{#1}}}}
\newcommand{\ten}[1] {\ensuremath{\boldsymbol{#1}}}
\newcommand{\tenfour}[1] {\ensuremath{\boldsymbol{\mathsf{#1}}}}
\newcommand{\comp}[1]{\begin{pmatrix}{ #1 }\end{pmatrix}}
\newcommand{\vv}{\ensuremath{\vec{v}}}
\newcommand{\vr}{\left(\ensuremath{\vec{r}}\right)}
\newcommand{\att}{\left(t\right)}
\newcommand{\tvv}{\ensuremath{\tilde{\vec{v}}}}
\newcommand{\ww}{\ensuremath{\vec{w}}}
\newcommand{\GG}{\ensuremath{\ten{G}}}
\newcommand{\FF}{\ensuremath{\boldsymbol{\mathfrak{F}}}}
\newcommand{\tEE}{\ensuremath{\tilde{\vec{E}}}}
\newcommand{\nn}{\ensuremath{\vec{n}}}
\renewcommand{\SS}{\ensuremath{\ten{S}}}
\newcommand{\tSS}{\ensuremath{\tilde{\ten{S}}}}
\newcommand{\cp}{\ensuremath{{c_p}}}
\newcommand{\cs}{\ensuremath{{c_s}}}
\newcommand{\tcs}{\ensuremath{{\tilde{c}_s}}}
\newcommand{\Vp}{\ensuremath{{V^e_p}}}
\newcommand{\Vs}{\ensuremath{{V^e_s}}}
\newcommand{\Ss}{\ensuremath{{S^e_s}}}
\newcommand{\nxnx}[1] {\ensuremath{\nn\times\left(\nn\times{#1}\right)}}

\newcommand{\dd}[2] {\ensuremath{\frac{\partial {#1}}{\partial {#2}}}}
\newcommand{\DD}[2] {\ensuremath{\frac{d {#1}}{d {#2}}}}
\newcommand{\Grad} {\ensuremath{\nabla}}
\newcommand{\Gradv}[1]{\ensuremath{\nabla_{#1}}}
\newcommand{\Div} {\ensuremath{\nabla\cdot}}
\newcommand{\Divv}[1] {\ensuremath{\nabla_{#1}\cdot}}
\newcommand{\Curl} {\ensuremath{\nabla\times}}
\newcommand{\Cten}{\ensuremath{\tenfour{C}}}

\newcommand{\eval}[2][\right]{\relax
  \ifx#1\right\relax \left.\fi#2#1\rvert}

\providecommand{\abs}[1]{\left\lvert#1\right\rvert}
\providecommand{\norm}[1]{\left\lVert#1\right\rVert}

\newcommand{\Nel} {\ensuremath{{N_\text{el}}}}
\newcommand{\Ned} {\ensuremath{{N_\text{ed}}}}
\newcommand{\De} {\ensuremath{{\mathsf{D}^e}}}
\newcommand{\Dep} {\ensuremath{{\mathsf{D}^{e'}}}}
\newcommand{\Dhat} {\ensuremath{\hat{\mathsf{D}}}}
\newcommand{\Dehat} {\ensuremath{{\hat{\mathsf{D}}^e}}}
\newcommand{\half} {\ensuremath{\frac{1}{2}}}
\newcommand{\IN} {\ensuremath{{\mathcal{I}_N}}}
\newcommand{\INe} {\ensuremath{{\mathcal{I}_{N_e}}}}

\newcommand{\mc}[1]{\mathcal{#1}}
\newcommand{\mcC}[1]{\mathcal{#1}}
\newcommand{\mb}[1]{\mathbf{#1}}
\newcommand{\mbb}[1]{\mathbb{#1}}
\newcommand{\mi}[1]{\mathit{#1}}
\newcommand{\mf}[1]{\mathfrak{#1}}
\newcommand{\Oi}[1]{\int_{\mc{D}} #1 \, d\vec{x}}
\newcommand{\TOi}[1]{\int_{0}^T \int_{\mc{D}} #1 \, d\vec{x}dt}
\newcommand{\Ti}[1]{\int_{0}^T #1 \, dt}
\newcommand{\TDei}[1]{\int_{0}^T\int_{\De} #1 \, d\vec{x}dt}
\newcommand{\TDeid}[1]{\int_{0}^T\int_{\De,N} #1 \, d\vec{x}dt}
\newcommand{\Dei}[1]{\int_{\De} #1 \, d\vec{x}}
\newcommand{\DeI}[1]{\int_{D} #1 \, d\vec{x}}
\newcommand{\DeiM}[1]{\int_{\Dhat} #1 \, d\vec{r}}
\newcommand{\DeMd}[1]{\int_{\De,N_e} #1 \, d\vec{x}}
\newcommand{\DeiMd}[1]{\int_{\Dhat,N_e} #1 \, d\vec{r}}
\newcommand{\DeMdN}[1]{\int_{D,N} #1 \, d\vec{x}}
\newcommand{\DeiMdN}[1]{\int_{\Dhat,N} #1 \, d\vec{r}}
\newcommand{\DeiMdNC}[2]{\int_{\Dhat,N_{#2}} #1 \, d\vec{r}}
\newcommand{\pDei}[1]{\int_{\partial \De} #1 \, d\vec{x}}
\newcommand{\pDeiM}[1]{\int_{\partial \Dhat} #1 \, d\vec{r}}
\newcommand{\pDeiMd}[1]{\int_{\partial \Dhat,N_e} #1 \, d\vec{r}}
\newcommand{\pDeiMdN}[1]{\int_{\partial \Dhat,N} #1 \, d\vec{r}}
\newcommand{\pDeiMdNC}[2]{\int_{\partial \Dhat,N_{#2}} #1 \, d\vec{r}}
\newcommand{\pMortar}[2]{\int_{\mc{M}_{#2}} #1 \, d\vec{r}}
\newcommand{\pMortard}[2]{\int_{\mc{M}^{#2},N_{#2}} #1 \, d\vec{r}}
\newcommand{\pMortardd}[3]{\int_{\mc{M}^{#2},N_{#3}} #1 \, d\vec{r}}
\newcommand{\pMortardh}[3]{\int_{\mc{M}_{#2},N_{#3}} #1 \, d\vec{r}}
\newcommand{\pMortardhn}[4]{\int_{{#2}_{#3},N_{#4}} #1 \, d\vec{r}}
\newcommand{\TpDei}[1]{\int_{0}^T\int_{\partial \De} #1 \, d\vec{x}}
\newcommand{\TpDeid}[1]{\int_{0}^T\int_{\partial \De,N} #1 \, d\vec{x}}
\newcommand{\Deid}[1]{\int_{\De,N_e} #1 \, d\vec{x}}
\newcommand{\DeidN}[1]{\int_{\De,N} #1 \, d\vec{x}}
\newcommand{\pDeidN}[1]{\int_{\partial \De,N} #1 \, d\vec{x}}
\newcommand{\pDeid}[1]{\int_{\partial \De,N_e} #1 \, d\vec{x}}
\newcommand{\nor}[1]{\left\| #1 \right\|}
\newcommand{\snor}[1]{\left| #1 \right|}
\newcommand{\LRp}[1]{\left( #1 \right)}
\newcommand{\LRs}[1]{\left[ #1 \right]}
\newcommand{\LRa}[1]{\left< #1 \right>}
\newcommand{\LRc}[1]{\left\{ #1 \right\}}
\newcommand{\pp}[2]{\frac{\partial #1}{\partial #2}}
\newcommand{\ppcp}[1]{\frac{\partial #1}{\partial \vec{c_p}}}
\newcommand{\ppcpj}[1]{\frac{\partial #1}{\partial c_p}}
\newcommand{\ppcpmf}[1]{\frac{\partial #1}{\partial \vec{\mf{c_p}}}}
\newcommand{\ppm}[1]{\frac{\partial #1}{\partial \vec{m}}}
\newcommand{\ppho}[3]{\frac{\partial^{#3} #1}{{\partial #2}^{#3}}}
\newcommand{\ppmi}[1]{\frac{\partial #1}{\partial m_i}}
\newcommand{\ppmj}[1]{\frac{\partial #1}{\partial m_j}}
\newcommand{\ppcpm}[1]{\frac{\partial #1}{\partial \vec{c_p}^-}}
\newcommand{\ppcpp}[1]{\frac{\partial #1}{\partial \vec{c_p}^+}}
\newcommand{\ppcs}[1]{\frac{\partial #1}{\partial \vec{c_s}}}
\newcommand{\ppcsm}[1]{\frac{\partial #1}{\partial \vec{c_s}^-}}
\newcommand{\ppcsp}[1]{\frac{\partial #1}{\partial \vec{c_s}^+}}
\newcommand{\td}[2]{\frac{d #1}{d #2}}
\newcommand{\Rvert}[1]{\left. #1 \right|}
\newcommand{\bs}[1]{\boldsymbol{#1}}
\newcommand{\ipDed}[3]{\LRp{ #1, #2}_{\De,N}^{#3}}

\newcommand{\bigO}{\mathcal{O}}

\newcommand{\iOm}[1]{\int_{\Omega} #1 \, d\Omega}
\newcommand{\iGb}[1]{\int_{\Gamma_{b}} #1 \, ds}
\newcommand{\iGs}[1]{\int_{\Gamma_{s}} #1 \, ds}
\newcommand{\iGsy}[1]{\int_{\Gamma_{s}} #1 \, ds(\mb{y})}
\newcommand{\RE}{\mbox{{I}}\hspace{-0.5ex}\mbox{{R}}}
\newcommand{\iC}[1]{\int_{0}^{2\pi} #1 \, d\theta}
\newcommand{\iGr}[1]{\int_{\Gamma_{\infty}} #1 \, ds}

\newcommand{\jump}[1] {\ensuremath{[\![{#1}]\!]}}
\newcommand{\diff}[1] {\ensuremath{\LRs{#1}}}
\newcommand{\average}[1] {\ensuremath{\LRc{\!\!\LRc{#1}\!\!}}}

\newcounter{tablecell}

\newcommand*{\numbercell}{%
  \stepcounter{tablecell}%
  \thetablecell. 
}

\newtheorem{propo}{Proposition}
\newtheorem{coro}{Corollary}
\newtheorem{theo}{Theorem}
\newtheorem{lem}{Lemma}

\newtheorem{defi}{definition}

\newtheorem{rema}{remark}

\newcommand{\figlab}[1]{\label{fig:#1}}
\newcommand{\eqnlab}[1]{\label{eq:#1}}
\newcommand{\theolab}[1]{\label{theo:#1}}
\newcommand{\corolab}[1]{\label{coro:#1}}
\newcommand{\propolab}[1]{\label{propo:#1}}
\newcommand{\lemlab}[1]{\label{lem:#1}}
\newcommand{\defilab}[1]{\label{defi:#1}}
\newcommand{\remalab}[1]{\label{rema:#1}}
\newcommand{\tablab}[1]{\label{tab:#1}}

\newcommand{\figref}[1]{\ref{fig:#1}}
\newcommand{\theoref}[1]{\ref{theo:#1}}
\newcommand{\defiref}[1]{\ref{defi:#1}}
\newcommand{\remaref}[1]{\ref{rema:#1}}
\newcommand{\cororef}[1]{\ref{coro:#1}}
\newcommand{\proporef}[1]{\ref{propo:#1}}
\newcommand{\lemref}[1]{\ref{lem:#1}}
\newcommand{\eqnref}[1]{\eqref{eq:#1}}
\newcommand{\alglab}[1]{\label{alg:#1}}
\newcommand{\algref}[1]{\ref{alg:#1}}
\newcommand{\seclab}[1]{\label{sect:#1}}
\newcommand{\secref}[1]{\ref{sect:#1}}
\newcommand{\tabref}[1]{\ref{tab:#1}}

\renewcommand{\u}{u}
\newcommand{\uk}{\u^k}
\newcommand{\ukp}{\u^{k+1}}
\renewcommand{\v}{v}
\newcommand{\vk}{{\v}^k}
\newcommand{\vkp}{{\v}^{k+1}}
\newcommand{\un}{\u_h}
\renewcommand{\e}{e}
\newcommand{\w}{w}
\renewcommand{\wb}{\mb{\w}}
\newcommand{\J}{J}
\newcommand{\Jb}{\mb{J}}
\newcommand{\q}{q}
\renewcommand{\r}{r}
\newcommand{\qh}{\hat{\q}}
\newcommand{\qs}{\q^*}
\newcommand{\qht}{\tilde{\qh}}
\newcommand{\qbar}{\overline{\q}}
\newcommand{\qb}{{\mb{\q}}}
\newcommand{\sig}{{\sigma}}
\newcommand{\thetah}{{\hat{\theta}}}
\newcommand{\thetas}{{\theta^*}}
\newcommand{\thetahn}{{\thetah_h}}
\newcommand{\sign}[1]{\text{ sgn}\!\LRp{#1}}
\newcommand{\sigh}{\hat{\sig}}
\newcommand{\sigs}{\sig^*}
\newcommand{\sigb}{{\bs{\sig}}}
\newcommand{\sigbn}{\sigb_h}
\newcommand{\sigbs}{\sigb^*}
\newcommand{\taub}{{\bs{\tau}}}
\newcommand{\kappab}{{\bs{\kappa}}}
\newcommand{\gamb}{{\bs{\gamma}}}
\newcommand{\ut}{\tilde{\u}}
\newcommand{\sigbt}{\tilde{\sigb}}
\newcommand{\vt}{\tilde{\v}}
\newcommand{\taubt}{\tilde{\taub}}
\newcommand{\ub}{\mb{\u}}
\newcommand{\ubp}{\mb{\u}^+}
\newcommand{\ubm}{\mb{\u}^-}
\newcommand{\vb}{\mb{\v}}
\newcommand{\um}{\u^-}
\newcommand{\up}{\u^+}

\newcommand{\m}{{m}}

\newcommand{\ud}{{\dot{\u}}}
\newcommand{\vd}{{\dot{\v}}}
\newcommand{\sigbd}{{\dot{\sigb}}}
\newcommand{\taubd}{\dot{\taub}}

\renewcommand{\d}{{d}}
\newcommand{\R}{{\mathbb{R}}}
\newcommand{\kap}{\kappa}
\newcommand{\F}{\mb{F}}
\newcommand{\f}{f}
\newcommand{\g}{g}
\newcommand{\fb}{\mb{\f}}
\newcommand{\gb}{\mb{\g}}
\newcommand{\gD}{g_D}
\newcommand{\pOmega}{\partial \Omega}
\newcommand{\K}{K}
\newcommand{\Kp}{K^+}
\newcommand{\Km}{K^-}
\newcommand{\pK}{{\partial \K}}
\newcommand{\pKin}{{\pK^{\text{in}}}}
\newcommand{\pKout}{{\pK^{\text{out}}}}
\newcommand{\Kj}{{\K_j}}
\newcommand{\Gh}{{\mc{E}_h}}
\newcommand{\Gho}{{\mc{E}_h^o}}
\newcommand{\Ghb}{{\mc{E}_h^{\partial}}}
\newcommand{\Poly}{{\mc{P}}}
\newcommand{\D}{{D}}
\newcommand{\p}{{p}}
\newcommand{\Ts}{{\mb{T}}}
\newcommand{\Vb}{{\mb{V}}}
\newcommand{\V}{{V}}
\newcommand{\Lam}{{\Lambda}}
\newcommand{\Lamb}{\bs{\Lambda}}
\newcommand{\mub}{\bs{\mu}}
\newcommand{\lambdab}{\bs{\lambda}}
\newcommand{\Th}{{\Ts_h}}
\newcommand{\Vh}{{\V_h}}
\newcommand{\Vbh}{{\Vb_h}}
\newcommand{\Lamh}{{\Lam_h}}
\newcommand{\Lambh}{{\Lamb_h}}
\newcommand{\Lamho}{{\overline{\Lam}_h}}
\newcommand{\Lambht}{{\Lamb_h^t}}
\newcommand{\Ltwo}{{L^2\LRp{\Omega}}}
\newcommand{\Ltwoe}{{L^2\LRp{\e}}}
\newcommand{\Ltw}{{L^2\LRp{\K}}}
\newcommand{\VhK}{{\V_h\LRp{\K}}}
\newcommand{\ThK}{{\Ts_h\LRp{\K}}}
\newcommand{\VbhK}{{\Vb_h\LRp{\K}}}
\newcommand{\Lamhe}{{\Lam_h\LRp{\e}}}
\newcommand{\Lambhe}{{\Lamb_h\LRp{\e}}}
\renewcommand{\L}{\mc{L}}
\renewcommand{\D}{\mc{D}}
\newcommand{\T}{\mc{T}}
\newcommand{\Tt}{\tilde{\T}}
\newcommand{\A}{\mb{A}}
\newcommand{\B}{\mb{B}}
\renewcommand{\D}{\mb{D}}
\newcommand{\Am}{\A^-}
\newcommand{\Ap}{\A^+}
\newcommand{\Aa}{\snor{\A}}
\newcommand{\Rb}{\mb{R}}
\newcommand{\C}{\mb{C}}
\renewcommand{\S}{\mb{S}}
\newcommand{\Sm}{\S^{-}}
\newcommand{\Sp}{\S^{+}}
\newcommand{\Spm}{\S^{\pm}}
\newcommand{\Q}{\mb{Q}}

\newcommand{\Lte}{{L^2\LRp{\Gh}}}
\newcommand{\n}{\mb{n}}
\newcommand{\nm}{\n^-}
\newcommand{\np}{\n^+}
\newcommand{\uh}{{\hat{\u}}}
\newcommand{\uhk}{{\uh}^k}
\newcommand{\uhkp}{{\uh}^{k+1}}
\newcommand{\us}{{\u^*}}
\newcommand{\uhn}{{\uh_h}}
\newcommand{\ubh}{{\hat{\ub}}}
\newcommand{\ubs}{{\ub^*}}
\newcommand{\ubht}{{\tilde{\ubh}}}
\newcommand{\uht}{{\tilde{\uh}}}
\newcommand{\uhtn}{{\uht_h}}
\newcommand{\uhd}{{\dot{\uh}}}
\newcommand{\sigbh}{{\hat{\sigb}}}
\newcommand{\Fh}{{\hat{\F}}}
\newcommand{\Fs}{\F^*}
\renewcommand{\c}{{c}}

\newcommand{\zb}{\mb{z}}
\renewcommand{\sb}{\mb{s}}
\newcommand{\rb}{\mb{r}}
\newcommand{\betab}{\bs{\beta}}
\newcommand{\veps}{{\varepsilon}}
\newcommand{\vepsb}{\bs{\veps}}
\newcommand{\Hb}{\mb{H}}
\newcommand{\Hbt}{\Hb^t}
\newcommand{\Eb}{\mb{E}}
\newcommand{\Ebt}{\Eb^t}
\newcommand{\hb}{\mb{h}}
\newcommand{\eb}{\mb{e}}
\newcommand{\Hbh}{\hat{\mb{H}}}
\newcommand{\Ebh}{\hat{\mb{E}}}
\newcommand{\Hbs}{{\mb{H}^*}}
\newcommand{\Ebs}{{\mb{E}^*}}
\newcommand{\Ebht}{\Ebh^t}
\newcommand{\Hbht}{\Hbh^t}
\newcommand{\Lb}{\mb{L}}
\newcommand{\Gb}{\mb{G}}
\newcommand{\Lbh}{\hat{\Lb}}
\newcommand{\Lbs}{\Lb^*}
\newcommand{\Ib}{\mb{I}}
\newcommand{\Sigb}{\bs{\Sigma}}
\newcommand{\Piu}{\Pi_\u}
\newcommand{\Pisigb}{\Pi_\sigb}
\newcommand{\Z}{Z}
\newcommand{\Y}{Y}
\newcommand{\rhou}{\rho\u}
\newcommand{\rhov}{\rho\v}
\newcommand{\E}{E}
\renewcommand{\P}{P}
\newcommand{\h}{H}

\newcommand{\mygreen}[1]{{\color[rgb]{0,.65,0} #1}}
\newcommand{\mywhite}[1]{{\color[rgb]{1.0,1.0,1.0} #1}}
\newcommand{\myred}[1]{{\color[rgb]{0.65,0.0,0.0} #1}}
\newcommand{\myblue}[1]{{\color[rgb]{0.0,0.0,1.0} #1}}
\newcommand{\mygray}[1]{{\color[rgb]{.15,.35,.60} #1}}
\newcommand{\mythemec}[1]{{\color[RGB]{51,108,121} #1}}
\newcommand{\mydarkred}[1]{{\color[rgb]{.6,.1,.1} #1}}
\newcommand{\mydarkblue}[1]{{\color[rgb]{.1,.1,.9} #1}}
\newcommand{\mygreenback}[1]{{\color[rgb]{.75,1,.75} #1}}
\newcommand{\myorange}[1]{{\color[rgb]{.76,.39,.13} #1}}
\newcommand{\mygrass}[1]{{\color[rgb]{.19,.64,.13} #1}}
\newcommand{\mysierp}[1]{{\color[RGB]{209,28,209} #1}}

\newcommand{\tanbui}[2]{\textcolor{blue}{#1} \textcolor{red}{#2}}

\newcommand{\vel}{\bs{\vartheta}}
\newcommand{\vels}{\vel^*}
\newcommand{\vele}{\vel^e}
\newcommand{\velh}{\hat{\vel}}
\newcommand{\velk}{\vel^k}
\newcommand{\velkp}{\vel^{k+1}}

\newcommand{\phis}{\phi^*}
\newcommand{\phie}{\phi^e}
\newcommand{\phih}{\hat{\phi}}
\newcommand{\phihk}{\phih^k}
\newcommand{\phik}{\phi^k}
\newcommand{\phikp}{\phi^{k+1}}

\maketitle

\begin{abstract}
We present a scalable and efficient iterative solver for high-order
hybridized discontinuous Galerkin (HDG) discretizations of hyperbolic
partial differential equations. It is an interplay between domain
decomposition methods and HDG discretizations. In particular, the method is a fixed-point
approach that requires only independent element-by-element local
solves in each iteration. As such, it is well-suited for current and
future computing systems with massive concurrencies. We rigorously
show that the proposed method is exponentially convergent in the number
of iterations for transport and linearized shallow water
equations. Furthermore, the convergence is independent of the solution
order. Various 2D and 3D numerical results for steady and
time-dependent problems are presented to verify our theoretical
findings.


\end{abstract}

\begin{keywords}
Iterative solvers, Discontinuous Galerkin methods, Hybridized Discontinuous Galerkin methods,
Shallow water equation, hyperbolic equations, Scalable solvers
\end{keywords}



\pagestyle{myheadings} \thispagestyle{plain}
\markboth{S. Muralikrishnan, M-B Tran and T. Bui-Thanh}{DDM-HDG}

\section{Introduction}

The discontinuous Galerkin (DG) method was originally developed by
Reed and Hill \cite{ReedHill73} for the neutron transport equation,
first analyzed in \cite{LeSaintRaviart74, JohnsonPitkaranta86}, and
then has been extended to other problems governed by partial
differential equations (PDEs) \cite{CockburnKarniadakisShu00}. Roughly
speaking, DG combines advantages of classical finite volume and finite
element methods. In particular, it has the ability to treat solutions
with large gradients including shocks, it provides the flexibility to
deal with complex geometries, and it is highly parallelizable due to
its compact stencil. As such, it has been adopted, for example, to solve large-scale
forward \cite{BreuerHeineckeRettenbergerEtAl13,
  WilcoxStadlerBursteddeEtAl10} and inverse
\cite{Bui-ThanhBursteddeGhattasEtAl12_gbfinalist} problems. 
However, for steady state problems or
time-dependent ones that require implicit time-integrators, DG methods
typically have many more (coupled) unknowns compared to the other
existing numerical methods, and hence more expensive in
general.

In order to mitigate the computational expense associated with DG
methods, Cockburn, coauthors, and others have introduced hybridizable
(also known as hybridized) discontinuous Galerkin (HDG) methods for
various types of PDEs including Poisson-type equation
\cite{CockburnGopalakrishnanLazarov09, CockburnGopalakrishnanSayas10,
  KirbySherwinCockburn12, NguyenPeraireCockburn09a,
  CockburnDongGuzmanEtAl09, EggerSchoberl10}, Stokes equation
\cite{CockburnGopalakrishnan09, NguyenPeraireCockburn10}, Euler and
Navier-Stokes equations, wave equations \cite{NguyenPeraireCockburn11,
  MoroNguyenPeraire11, NguyenPeraireCockburn11b,
  LiLanteriPerrrussel13, NguyenPeraireCockburn11a, GriesmaierMonk11,
  CuiZhang14}, to name a few. The upwind HDG framework proposed in
\cite{Bui-Thanh15, Bui-Thanh15a, Bui-Thanh15b} provides a unified and
a systematic construction of HDG methods for a large class of PDEs. In
HDG discretizations, the coupled unknowns are single-valued traces
introduced on the mesh skeleton, i.e. the faces, and the resulting
matrix is substantially smaller and sparser compared to standard DG
approaches. Once they are solved for, the usual DG unknowns can be
recovered in an element-by-element fashion, completely independent of
each other.  Nevertheless, the trace system is still a bottleneck for
practically large-scale applications, where complex and high-fidelity
simulations involving features with a large range of spatial and
temporal scales are necessary.

Meanwhile, Schwarz-type domain decomposition methods (DDMs) have been
introduced as procedures to parallelize and solve partial differential
equations numerically, where each iteration involves the solutions of
the original equations on smaller subdomains
\cite{Lions:1987:OSA,Lions:1989:OSA,Lions:1990:OSA}. Among the many
DDMs, Schwarz waveform relaxation methods and optimized Schwarz
methods
\cite{Halpern:OSM:2009,BennequinGanderHalpern:AHB:2009,HalpernSzeftel:2009:NOS,GanderGouarinHalpern:2011:OSW,GanderHajian:2015:ASM,Binh1},
have attracted substantial attention over the past decade since they
can be adapted to the physics of the underlying problems and thus lead
to very efficient parallel solvers for challenging problems. We view
the HDG method as an extreme DDM approach in which each subdomain is
an element.

While either HDG community or DDM community can contribute
individually towards advancing its own field, the potential for a true
breakthrough may lie in bringing together the advances from both sides
and in exploiting opportunities at their interfaces.  In this paper,
we blend the HDG method and optimized Schwarz idea to produce a
efficient and scalable iterative approach for HDG methods. One of the
main features of the proposed approach is that it has exponential
convergence rate, and for that reason we term it as eHDG. The method
can be viewed as a fixed-point approach that requires only independent
element-by-element local solves in each iteration. As such, it is
well-suited for current and future computing systems with massive
concurrencies. We rigorously show that the proposed method is
exponentially convergent in the number of iterations for transport and
linearized shallow water equations. Furthermore, the convergence is
independent of the solution order. The theoretical findings will be
verified on various 2D and 3D numerical results for steady and
time-dependent problems.

Let us mention that in \cite{GanderHajian:2015:ASM}, Schwarz methods for the hybridizable interior penalty (IPH) method have also been introduced. The methods have been proposed entirely at
the discrete level and thus holds for arbitrary interfaces between two subdomains. It is proved that for an arbitrary two-subdomain decomposition the Schwarz algorithms  have a convergence factor $1-O(h)$, and $1 -O(\sqrt{h}),$ which means the algorithms converge slower and slower when we refine the mesh.

\section{Notations for HDG discretizations}
\seclab{HDG} In this section we introduce common notations and conventions to be
used in the following sections where we propose and rigorously analyze
the eHDG approach for scalar and systems of hyperbolic PDEs in both steady
and time-dependent cases.  Let us partition an open and bounded domain $\Omega \in \R^\d$ into
$\Nel$ non-overlapping elements $\Kj, j = 1,\hdots,\Nel$ with
Lipschitz boundaries such that $\Omega_h := \cup_{j=1}^\Nel \Kj$ and
$\overline{\Omega} = \overline{\Omega}_h$. Here, $h$ is defined as $h
:= \max_{j\in \LRc{1,\hdots,\Nel}}diam\LRp{\Kj}$. We denote the
skeleton of the mesh by $\Gh := \cup_{j=1}^\Nel \partial K_j$,
the set of all (uniquely defined) faces $\e$. We conventionally identify $\nm$ as the normal
vector on the boundary $\pK$ of element $\K$ (also denoted as $\Km$)
and $\np = -\nm$ as the normal vector of the boundary of a neighboring
element (also denoted as $\Kp$). Furthermore, we use $\n$ to denote
either $\nm$ or $\np$ in an expression that is valid for both cases,
and this convention is also used for other quantities (restricted) on
a face $\e \in \Gh$. 

For simplicity in writing we define $\LRp{\cdot,\cdot}_\K$ as the
$L^2$-inner product on a domain $\K \in \R^\d$ and
$\LRa{\cdot,\cdot}_\K$ as the $L^2$-inner product on a domain $\K$ if
$\K \in \R^{\d-1}$. We shall use $\nor{\cdot}_{\K} :=
\nor{\cdot}_{\Ltw}$ as the induced norm for both cases and the
particular value of $\K$ in a context will indicate which inner
product the norm is coming from. We also denote the $\veps$-weighted
norm of a function $\u$ as $\nor{\u}_{\veps, \K} :=
\nor{\sqrt{\veps}\u}_{\K}$ for any positive $\veps$.  We shall use
boldface lowercase letters for vector-valued functions and in that
case the inner product is defined as $\LRp{\ub,\vb}_\K :=
\sum_{i=1}^m\LRp{\ub_i,\vb_i}_\K$, and similarly $\LRa{\ub,\vb}_\K :=
\sum_{i = 1}^m\LRa{\ub_i,\vb_i}_\K$, where $\m$ is the number of
components ($\ub_i, i=1,\hdots,\m$) of $\ub$.  Moreover, we define
$\LRp{\ub,\vb}_\Omega := \sum_{\K\in \Omega_h}\LRp{\ub,\vb}_\K$ and
$\LRa{\ub,\vb}_\Gh := \sum_{\e\in \Gh}\LRa{\ub,\vb}_\e$ whose
induced (weighted) norms are clear, and hence their definitions are
omitted. We  employ boldface uppercase letters, e.g. $\mb{L}$, to
denote matrices and tensors. In addition, subscripts
are used to denote the components of vectors, matrices, and tensors.

We define $\Poly^\p\LRp{\K}$ as the space of polynomials of degree at
most $\p$ on a domain $\K$. Next, we introduce two discontinuous
piecewise polynomial spaces
\begin{align*}
\Vbh\LRp{\Omega_h} &:= \LRc{\vb \in \LRs{L^2\LRp{\Omega_h}}^2:
  \eval{\vb}_{\K} \in \LRs{\Poly^\p\LRp{\K}}^2, \forall \K \in \Omega_h}, \\
\Lambh\LRp{\Gh} &:= \LRc{\lambdab \in \LRs{\Lte}^2:
  \eval{\lambdab}_{\e} \in \LRs{\Poly^\p\LRp{\e}}^2, \forall \e \in \Gh},
\end{align*}
and similar spaces for $\VbhK$ and $\Lambhe$ by replacing $\Omega_h$ with
$\K$ and $\Gh$ with $\e$. For scalar-valued functions,
we denote the corresponding spaces as
\begin{align*}
\Vh\LRp{\Omega_h} &:= \LRc{\v \in L^2\LRp{\Omega_h}:
  \eval{\v}_{\K} \in \Poly^\p\LRp{\K}, \forall \K \in \Omega_h}, \\
\Lamh\LRp{\Gh} &:= \LRc{\lambda \in \Lte:
  \eval{\lambda}_{\e} \in \Poly^\p\LRp{\e}, \forall \e \in \Gh}.
\end{align*}

\section{Construction of eHDG methods for linear hyperbolic PDEs}
In this section, we define eHDG methods for scalar and system of
hyperbolic PDEs. For the clarity in exposition, we consider the
transport equation and linearized shallow water system, and
extension of the proposed approach to other hyperbolic PDEs is straightforward. 
To begin, let us consider the transport equation
\begin{subequations}
\eqnlab{transport}
\begin{align}
\betab \cdot \Grad \u &= f \quad \text{ in },
\Omega, \\
\u &= \g \quad \text{ on } \pOmega^-,
\end{align}
\end{subequations}
where $\pOmega^-$ is the inflow part of the boundary $\pOmega$.  An
upwind HDG discretization \cite{Bui-Thanh15} for
\eqnref{transport} consists of the following local equation for each element $\K$
\begin{equation}
\eqnlab{transportLocal}
-\LRp{\u,\Div\LRp{\betab \v}}_\K  + \LRa{\betab\cdot\n\u + \snor{\betab\cdot \n}\LRp{\u -\uh},\v}_\pK= \LRp{f,\v}_\K, \quad \forall \v \in \Vh\LRp{\K},
\end{equation}
and conservation conditions on all edges $\e$ in the mesh skeleton $\Gh$:
\[
\LRa{\jump{\betab\cdot\n\u + \snor{\betab\cdot \n}\LRp{\u -\uh}},\mu}_\e = 0, \quad \forall \mu \in \Lamh\LRp{\e}.
\]



Inspired by the upwind HDG approach \cite{Bui-Thanh15} and the
optimized Schwarz method \cite{Binh1}, we introduce an eHDG iterative
method for the transport equation \eqnref{transport} as in Algorithm
\ref{al:DDMhyperbolic}. In particular, the approximate solution
$\u^{k+1}$ at the $(k+1)$th iteration restricted on element $\K$ is
defined as the solution of the following local equation, $\forall \v
\in \Vh\LRp{\K}$,
\begin{equation}
\eqnlab{transportLocalk1}
- \LRp{\ukp,\Div\LRp{\betab \v}}_\K  + \LRa{\betab\cdot\n\ukp + \snor{\betab\cdot \n}\LRp{\ukp -\uhk},\v}_\pK= \LRp{f,\v}_\K,
\end{equation}
where, by introducing the average operator as $2\average{\LRp{\cdot}} :=
\LRp{\cdot}^- + \LRp{\cdot}^+$, we define
\begin{equation}
\eqnlab{transportLocale}
\uhk := \average{\uk\sign{\betab\cdot\n}} + \average{\uk}.
\end{equation}

\begin{algorithm}
  \begin{algorithmic}[1]
    \ENSURE Given initial guess $\u^0$, compute the initial
    trace $\uh^0$ using \eqnref{transportLocale}.
      \WHILE{not converged} 
      \STATE Solve the local equation \eqnref{transportLocalk1} for $\ukp$ using trace $\uhk$
      \STATE Compute $\uh^{k+1}$ using \eqnref{transportLocale}.
      \STATE Check convergence. If yes, {\bf exit}, otherwise {\bf continue}
    \ENDWHILE
  \end{algorithmic}
  \caption{eHDG solver for transport equation \eqnref{transport}}
  \label{al:DDMhyperbolic}
\end{algorithm}

Since \eqnref{transport} is linear, it is sufficient to show that eHDG converges for the homogeneous equation with zero forcing $\f$ and zero boundary condition $\g$.
\begin{theorem}\theolab{DDMConvergence} Assume $-\Div \betab \ge \alpha > 0$, i.e. \eqnref{transport} is well-posed. The above eHDG for homogeneous transport equation \eqnref{transport} converges exponentially with respect to the number of iterations $k$. In particular, there exist $J \le \Nel$ such that 
\begin{align}
\label{ThmExponentialDecay}
\nor{\uk}^2_{\frac{-\Div \betab}{2}, \Omega_h} +
\nor{\uk}_{\snor{\betab\cdot\n},\Gh}^2&\le
\frac{C(k)}{2^k}\nor{\u^0}_{\snor{\betab\cdot\n},\Gh}^2,
\end{align}
 where $C(k)$ is a polynomial in $k$ of order at most $J$ and is independent of $\p$.
\end{theorem}


For time-dependent transport equation, we discretize the spatial
operator using HDG and time using backward Euler method (for
simplicity). The eHDG in this case is almost identical to the one for
steady state equation except that we now have an additional $L^2$-term
in the local equation \eqnref{transportLocalk1}. 

We next consider the following oceanic linearized shallow water systems \cite{GiraldoWarburton08}
\begin{equation}
    \pp{}{t}\LRp{
      \begin{array}{c}
        \phi \\
       \Phi \u \\
        \Phi \v
      \end{array}
    } + 
    \pp{}{x}\LRp{
      \begin{array}{c}
        \Phi  \u \\
        \Phi\phi \\
        0
      \end{array}
    } + 
    \pp{}{y}\LRp{
      \begin{array}{c}
        \Phi \v \\
        0 \\
        \Phi\phi
      \end{array}
    } = 
    \LRp{
      \begin{array}{c}
        0 \\
        f\Phi\v - \gamma \Phi\u + \frac{\tau_x}{\rho} \\
       -f\Phi\u - \gamma \Phi\v + \frac{\tau_y}{\rho}
      \end{array}
    }
\eqnlab{linearizedShallow}
\end{equation}
where $\phi = g \h$ is the geopotential height with $g$ and $\h$ being
the gravitational constant and the perturbation of the free surface
height, $\Phi > 0$ is a constant mean flow geopotential height, $\vel := \LRp{\u,\v}$ is
the perturbed velocity, $\gamma \ge 0$ is the bottom friction, 
$\bs{\tau}:=\LRp{\tau_x,\tau_y}$ is the wind stress, and $\rho$ is the density of
the water. Here, $f = f_0 + \beta \LRp{y - y_m}$ is the
Coriolis parameter,  where $f_0$, $\beta$, and $y_m$ are given constants.

Again, for simplicity of the exposition and analysis, let us employ
the backward Euler discretization for temporal derivatives and HDG for spatial ones.
Since the unknowns of interest
are those at the $(m+1)$th time step, we can suppress the time index
for clarity of the exposition.  Furthermore, since the system
\eqnref{linearizedShallow} is linear, a similar argument as in the previous
sections shows that it is sufficient to consider homogeneous system with zero initial condition, boundary condition, and forcing. An eHDG algorithm can be proposed
for the  homegeneous system as follows
\begin{subequations}
\eqnlab{localSolverDD}
\begin{align*}
&\LRp{\frac{\phikp}{\Delta t},\varphi_1}_K - \LRp{\Phi\velkp,\Grad\varphi_1}_K + \LRa{\Phi \velkp \cdot \n + \sqrt{\Phi}\LRp{\phikp -\phihk},\varphi_1}_\pK = 0, \\
&\LRp{\frac{{\Phi\ukp}}{\Delta t},\varphi_2}_K - \LRp{\Phi\phikp,\pp{\varphi_2}{x}}_K + \LRa{\Phi\phihk\n_1,\varphi_2}_\pK = \LRp{f\Phi\vkp - \gamma \Phi\ukp, \varphi_2}_K, \\
&\LRp{\frac{{\Phi\vkp}}{\Delta t},\varphi_3}_K - \LRp{\Phi\phikp,\pp{\varphi_3}{y}}_K + \LRa{\Phi\phihk\n_2,\varphi_3}_\pK = \LRp{-f\Phi\ukp - \gamma \Phi\vkp , \varphi_3}_K,
\end{align*}
\end{subequations}
where $\varphi_1, \varphi_2$ and $\varphi_3$ are the test functions, and similar to the transport equation we define
\begin{equation*}
\phihk= \average{\phik} + \sqrt{\Phi}\average{\velk\cdot\n}.
\end{equation*}
Our goal is to show that $\LRp{\phikp,\Phi\velkp}$ converges to
zero. To that end, let us define
\begin{eqnarray}\eqnlab{ContractionConstant}
\mathcal{C}:=\frac{\mathcal{A}}{\mathcal{B}}, \quad \mathcal{A}
:=\max\left\{\frac{\Phi+\sqrt\Phi}{2},\frac{(1+\sqrt\Phi)}{2}\right\},
\end{eqnarray}
and
\begin{eqnarray*}
\mathcal{B}
:=\min\left\{\left(\frac{h}{\Delta t(\p+1)(\p+2)}+\frac{\sqrt\Phi -\Phi}{2}\right),\left(\frac{ h}{\Delta t(\p+1)(\p+2)}+\frac{(2\gamma-1-\sqrt\Phi)}{2}\right)\right\}.
\end{eqnarray*}
We also need the following norms:
\begin{align}
\nor{\LRp{\phik,\velk}}^2_{\Omega_h}&:=\nor{\phik}_{\Omega_h}^2 + \nor{\velk}_{\Phi, \Omega_h}^2, \quad \nor{\LRp{\phik,\velk}}^2_\Gh&:=\nor{\phik}_{\Gh}^2 + \nor{\velk}_{\Phi, \Gh}^2.\nonumber
\end{align}

\begin{theorem}\label{ThmConvergenceShallowWater} If the mesh size $h$, the time step $\Delta t$ and the order $\p$ are chosen such that $\mc{B} > 0$ and $\mathcal{C}<1$,
then the approximate solution at the $k$th iteration
$\LRp{\phik,\velk}$ decays exponentially in the following sense
\begin{align}\nonumber
\nor{\LRp{\phik,\velk}}^2_\Gh &\leq \mathcal{C}^k \nor{\LRp{\phi^0,\vel^0}}^2_\Gh, \\
\nor{\LRp{\phik,\velk}}^2_{\Omega_h} &\leq \Delta t \mc{A}\LRp{\mc{C}+1}\mc{C}^k\nor{\LRp{\phi^0,\vel^0}}^2_\Gh, \nonumber
\end{align}
where $\mc{C}$ is defined in \eqnref{ContractionConstant}.
\end{theorem}

\section{Numerical results}
In this section various numerical results supporting the theoretical results are provided for
2D and 3D transport equations and the linearized shallow water equation.

\subsection{2D steady state transport equation with smooth solution}
\label{ex1_2dsteady}
In this example we choose $\betab=(y,x)$. Also we take the forcing and the exact solution 
to be of the following form:
\begin{subequations}
\eqnlab{exact2dsteady}
\begin{align}
	u^e&=\frac{1}{\pi}\sin(\pi x)\cos(\pi y), \\
	f&=y\cos(\pi x)\cos(\pi y)-x\sin(\pi x)\sin(\pi y).
\end{align}
\end{subequations}
Here the domain $\Omega$ is $[0,1]\times[0,1]$ with 
$x=0$ and $y=0$ as inflow boundaries. A structured quadrilateral mesh is used for all the 
numerical simulations performed.
\begin{figure}[h!t]
\vspace{-5mm}
\subfigure[2D steady transport]{
\includegraphics[trim=3.5cm 7.75cm 6.75cm 9.25cm,clip=true,width=0.3\columnwidth]{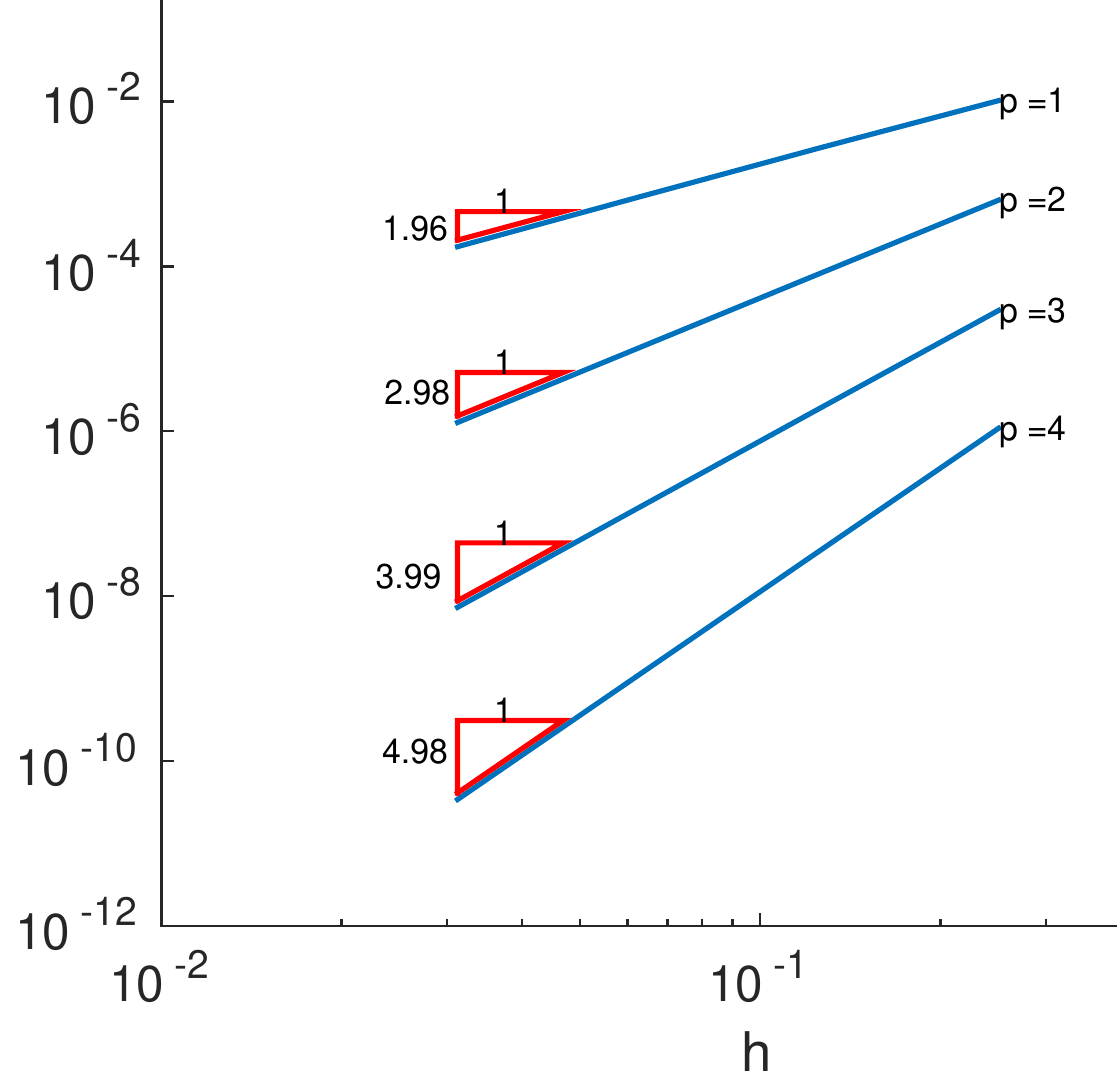}
}
\subfigure[3D steady transport]{
\includegraphics[trim=3.5cm 8cm 5.75cm 10cm,clip=true,width=0.3\columnwidth]{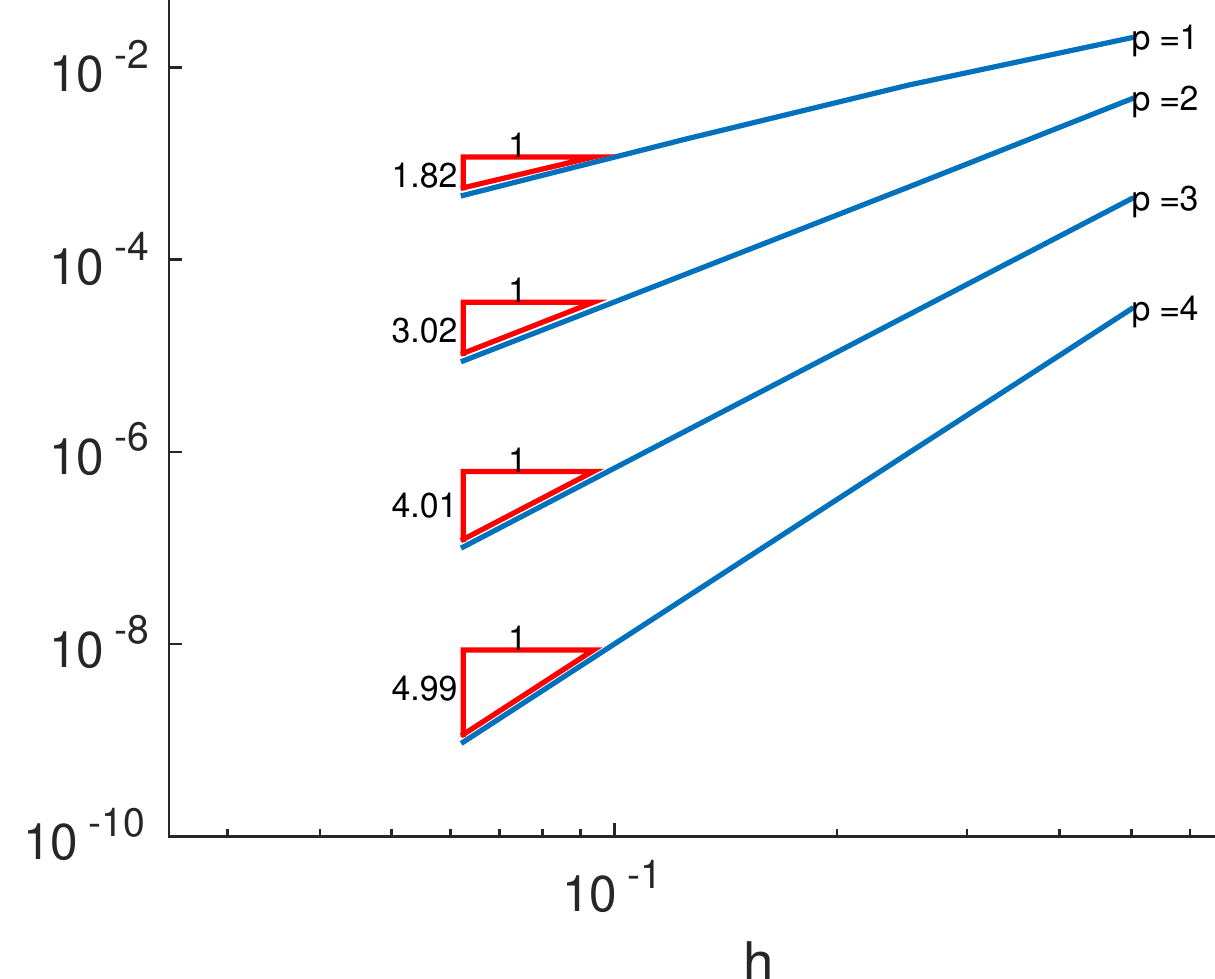}
\label{l2err_comb-b}
}
\subfigure[shallow water equation]{
\includegraphics[trim=3.5cm 7.75cm 5.5cm 9.25cm,clip=true,width=0.3\columnwidth]{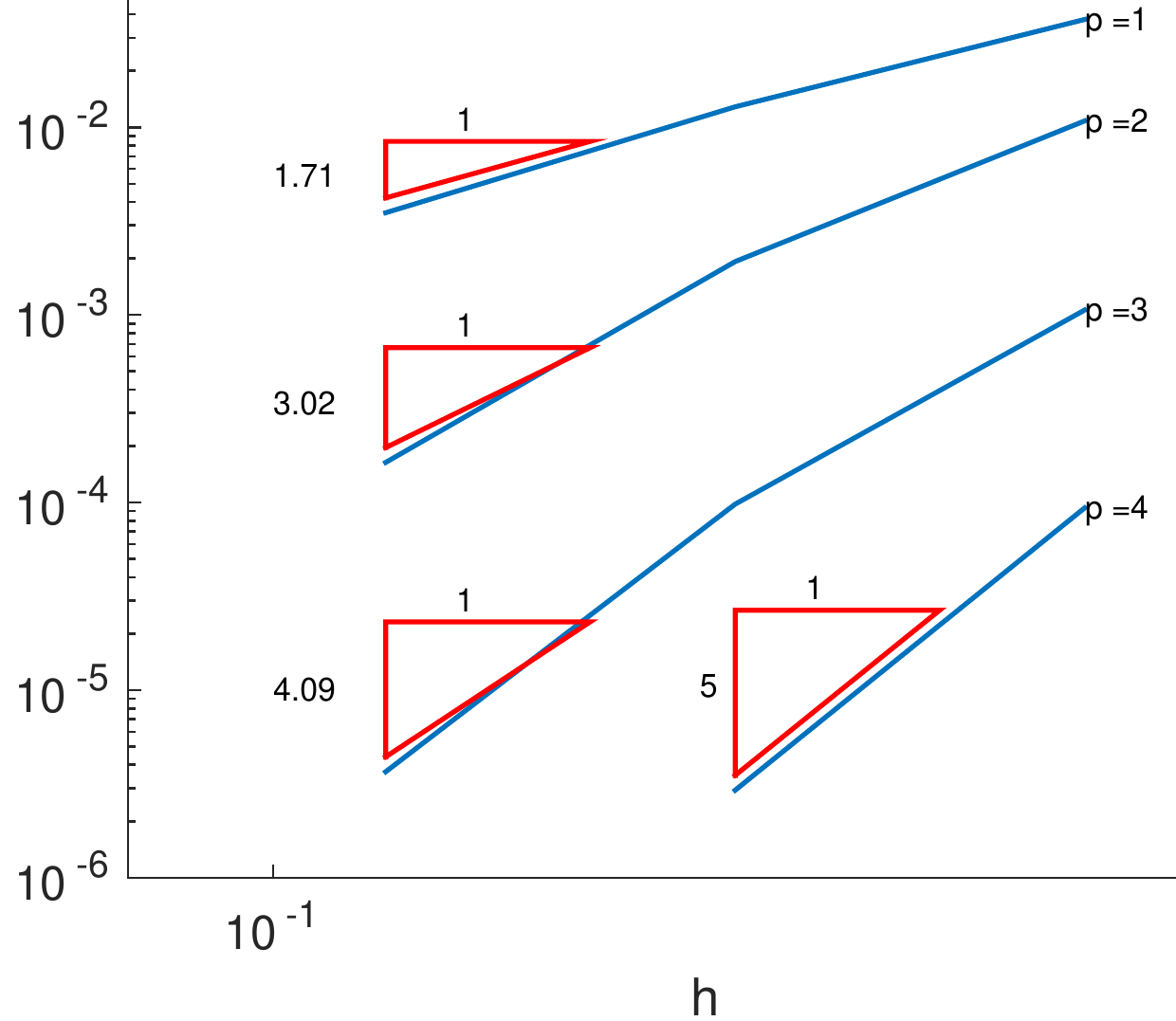}
\label{l2err_comb-c}
}
\caption{ h-convergence of the HDG method using the eHDG solver.} 
\label{l2err_comb}
\end{figure}

\begin{figure}[h!t!b!]
\vspace{-7mm}
\subfigure[Error history for p=1,2]{
\includegraphics[trim=2cm 7cm 2.5cm 7.25cm,clip=true,width=0.5\textwidth]{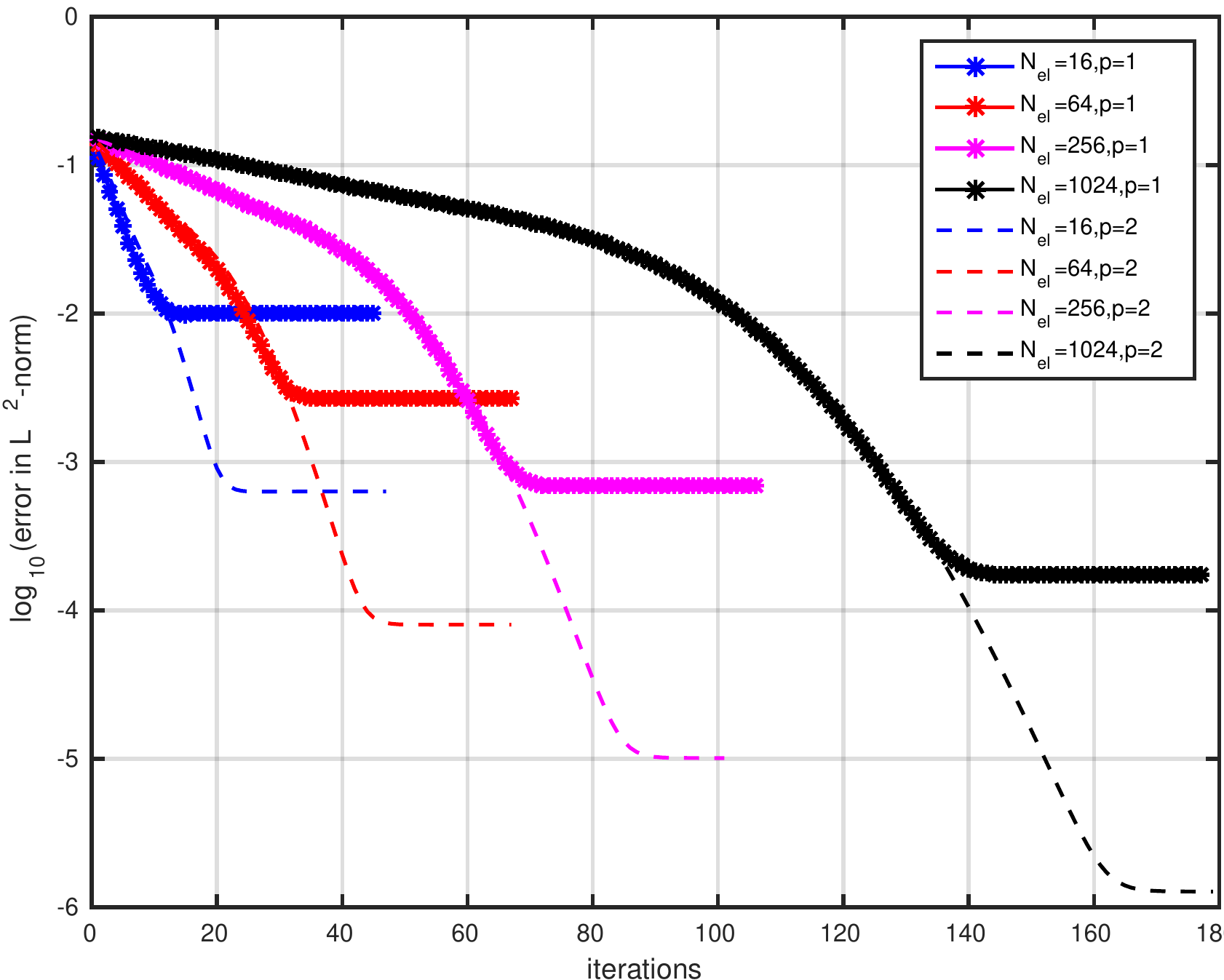}
}
\subfigure[Error history for p=3,4]{
\includegraphics[trim=1.75cm 7cm 2.5cm 7.25cm,clip=true,width=0.5\textwidth]{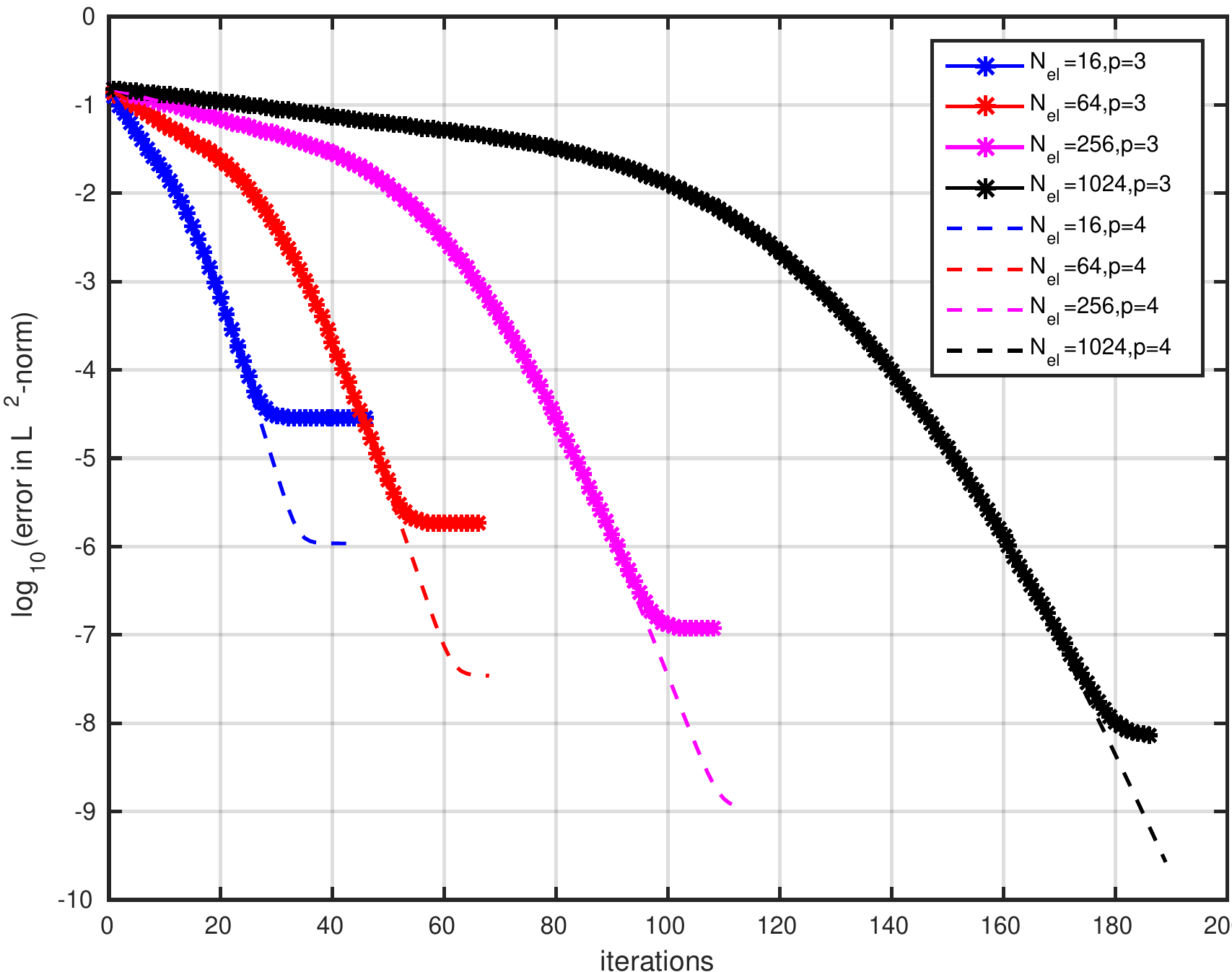}
}
\caption{Convergence of eHDG for different $h$ and $p$ for 2D transport smooth solution.} 
\label{allp_2dsteady}
\end{figure}
\vspace{-6mm}
\begin{minipage}{\linewidth}
\centering
\captionof{table}{\it{Iterations for 2D transport equation with smooth and discontinuous solutions and 3D steady state transport equation}} 
\label{tab:Comparison_2dall_3d}
\begin{tabular}{ | r || r || r || c | c | c | }
\hline
$N_{el}(2D)$ & $N_{el}(3D)$ & p & 2D smooth & 2D discontinuous & 3D steady \\
\hline
16 & 8 & 1 & 45 & 59 & 33 \\
\hline
64 & 64 & 1 & 67 & 84 & 51 \\
\hline
256 & 512 & 1 & 107 & 129 & 79 \\
\hline
1024 & 4096 & 1 & 177 & 209 & 130 \\
\hline
\hline
16 & 8 & 2 & 47 & 61 & 39 \\
\hline
64 & 64 & 2 & 67 & 87 & 51 \\
\hline
256 & 512 & 2 & 101 & 133 & 76 \\
\hline
1024 & 4096 & 2 & 179 & 214 & 131 \\
\hline
\hline
16 & 8 & 3 & 46 & 65 & 39 \\
\hline
64 & 64 & 3 & 66 & 92 & 49 \\
\hline
256 & 512 & 3 & 108 & 135 & 79 \\
\hline
1024 & 4096 & 3 & 186 & 211 & 136 \\
\hline
\hline
16 & 8 & 4 & 45 & 66 & 35 \\
\hline
64 & 64 & 4 & 68 & 90 & 51 \\
\hline
256 & 512 & 4 & 112 & 128 & 83 \\
\hline
1024 & 4096 & 4 & 189 & 198 & 143 \\
\hline
\end{tabular}
\end{minipage}

Figure \ref{l2err_comb} shows the $h$-convergence of the HDG discretization with eHDG iterative solver. The convergence is optimal i.e. $(p+1)$ for a polynomial order $p$. The tolerance criteria for the eHDG solver is set as follows:
\begin{equation}
\eqnlab{tolerancecriteria}
|\norm{u^k-u_e}_{\Ltwo}-\norm{u^{k-1}-u_e}_{\Ltwo}|<10^{-10}.
\end{equation}
Thus the succesive difference in $L^2$ norm of error between numerical solution and exact solution is used as a criteria 
for tolerance in this case.

Figure \ref{allp_2dsteady} shows the convergence history of the eHDG
solver in the log-linear scale. As proved in Theorem
\theoref{DDMConvergence} the eHDG is exponential convergent in the
iteration $k$. Also the stagnation region observed near the end of
each curve is due to the fact that for a particular mesh size $h$ and
polynomial order $p$ we can achieve only as much accuracy as
prescribed by the HDG discretization error and cannot go beyond
that. The numerical results for different solution orders also verify
the fact that the convergence of eHDG method is independent of the
polynomial order $p$. This can also be seen from the $4$th column of
Table \ref{tab:Comparison_2dall_3d}.
\subsection{2D steady state transport equation with discontinuous solution}
In this case we take  $f=0$ and $\betab=(1+sin(\pi y/2),2)$. The domain $\Omega$ is $[0,2]\times[0,2]$ and the inflow boundary condition is given as 
\[
g = 
\left\{
\begin{array}{ll}
1 & x = 0, 0 \le y \le 2 \\
\sin^6\LRp{\pi x} & 0< x \le 1, y = 0 \\
0 & 1 \le x \le 2, y = 0
\end{array}
\right.
.
\]
We choose a slight different stopping criteria to avoid the exact solution:
\[
\norm{u^k-u^{k-1}}_{\Ltwo}<10^{-10}.
\]
The evolution of solution with iterations obtained for $32\times32$
elements and polynomial order $4$ is shown in Figure
\figref{2ddiscontiterfig1}. As shown from the $5$th column of Table \ref{tab:Comparison_2dall_3d}, due to the discontinuity, the eHDG solver
takes a slightly more iterations compared to the smooth solution case,
but the number of iteration is still (almost) independent of the
solution order. Also we observe that the solution evolves from inflow
to outflow. This can be proved rigorously, but for
the space limitation, the proof is omitted.
\begin{figure}[h!t!b!]
  \subfigure[$\u$ at $iteration = 16$]{
    \includegraphics[trim=2.5cm 1.5cm 2.5cm 1cm,clip=true,width=0.3\columnwidth]{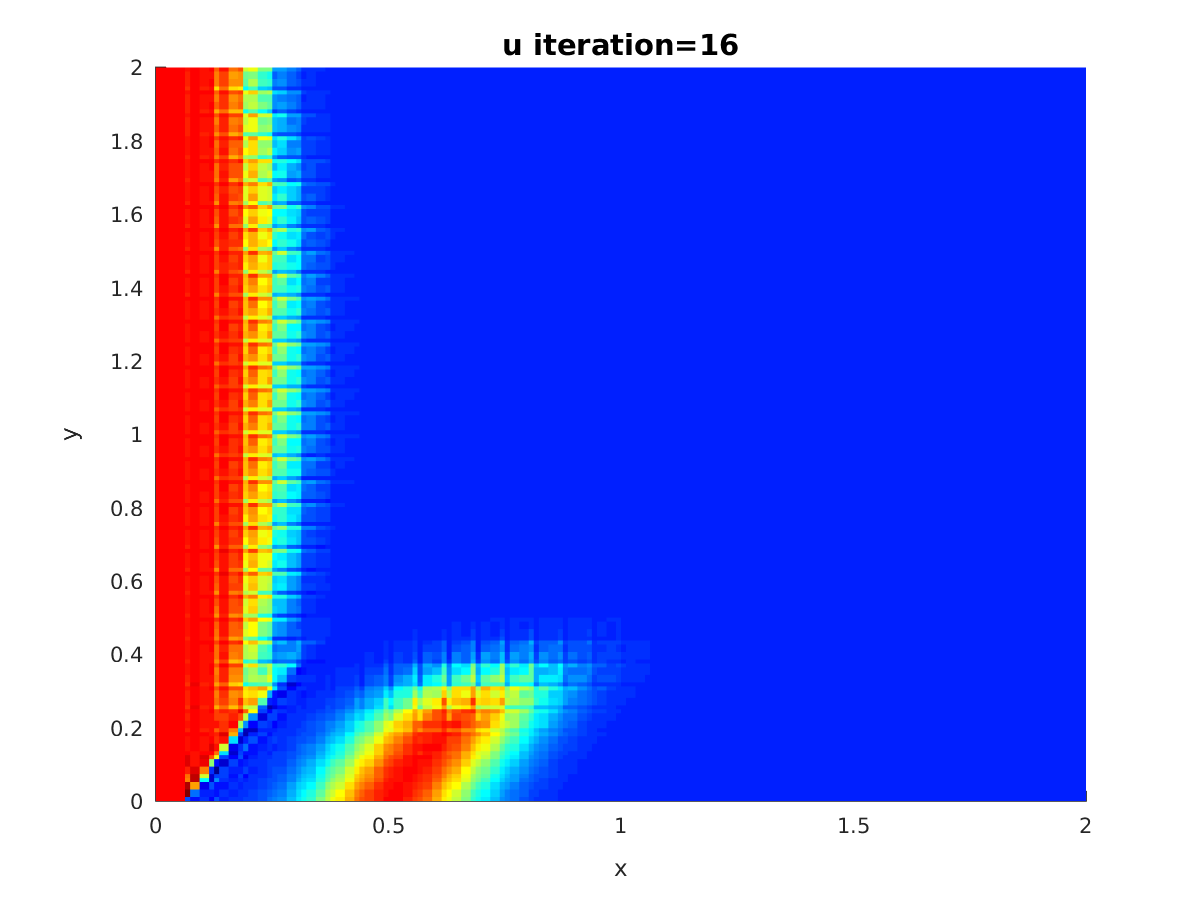}
  }
  \subfigure[$\u$ at $iteration = 64$]{
    \includegraphics[trim=2.5cm 1.5cm 2.5cm 1cm,clip=true,width=0.3\columnwidth]{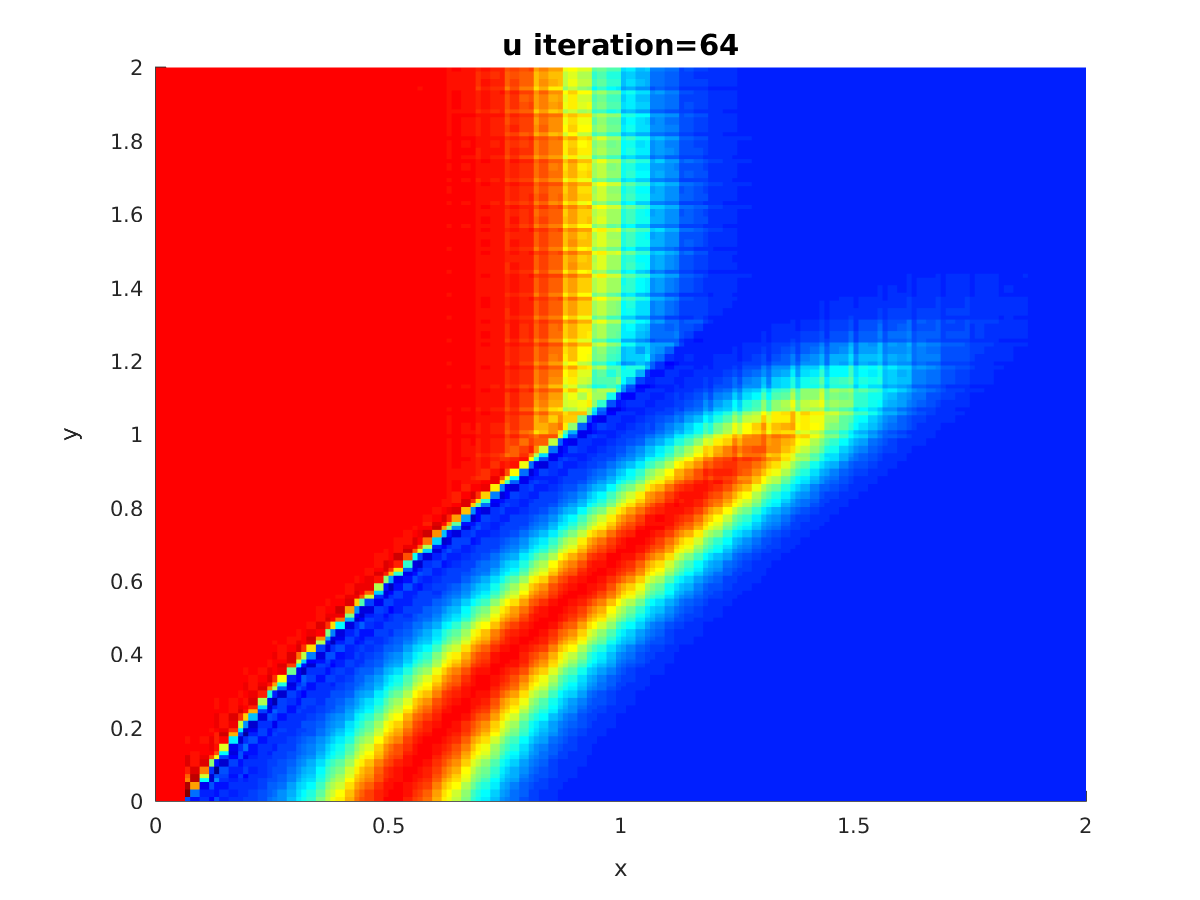}
  }
  \subfigure[$\u$ at $iteration = 192$]{
    \includegraphics[trim=2.5cm 1.5cm 2.5cm 1cm,clip=true,width=0.3\columnwidth]{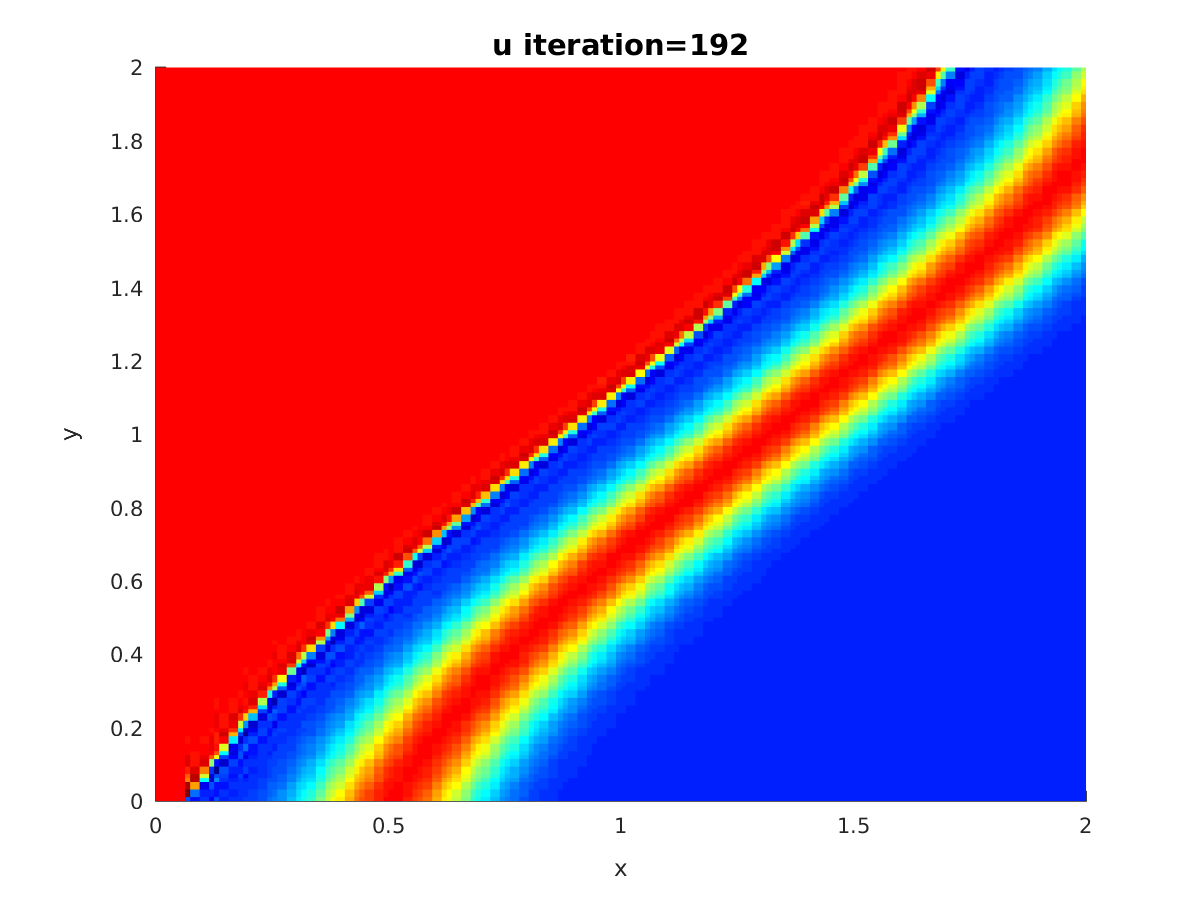}
  }
 \caption{Evolution of solution with respect to iterations for upwind HDG}
  \figlab{2ddiscontiterfig1}
\end{figure}

\subsection{3D steady state transport equation}

In this example we choose $\betab=(z,x,y)$ in \eqnref{transport} and  the following exact solution:
$$u^e=\frac{1}{\pi}\sin(\pi x)\cos(\pi y)\sin(\pi z).$$
The forcing is selected in such a way that it corresponds to the exact solution selected.
Here, the domain is $[0,1]\times[0,1]\times[0,1]$ with 
$x=0$, $y=0$ and $z=0$ as inflow boundaries. A structured $16\times16\times16$  hexahedral mesh is used for all
simulations. The tolerance critera used is same
as in section \ref{ex1_2dsteady}.
Similar to the 2D example in section \ref{ex1_2dsteady}, we obtain the optimal convergence
rates as shown in Figure \ref{l2err_comb-b}. 
\begin{figure}[h!t!b!]
\vspace{-6mm}
\subfigure[Error history for p=1,2]{
\includegraphics[trim=1.75cm 7cm 2.5cm 7.25cm,clip=true,width=0.5\textwidth]{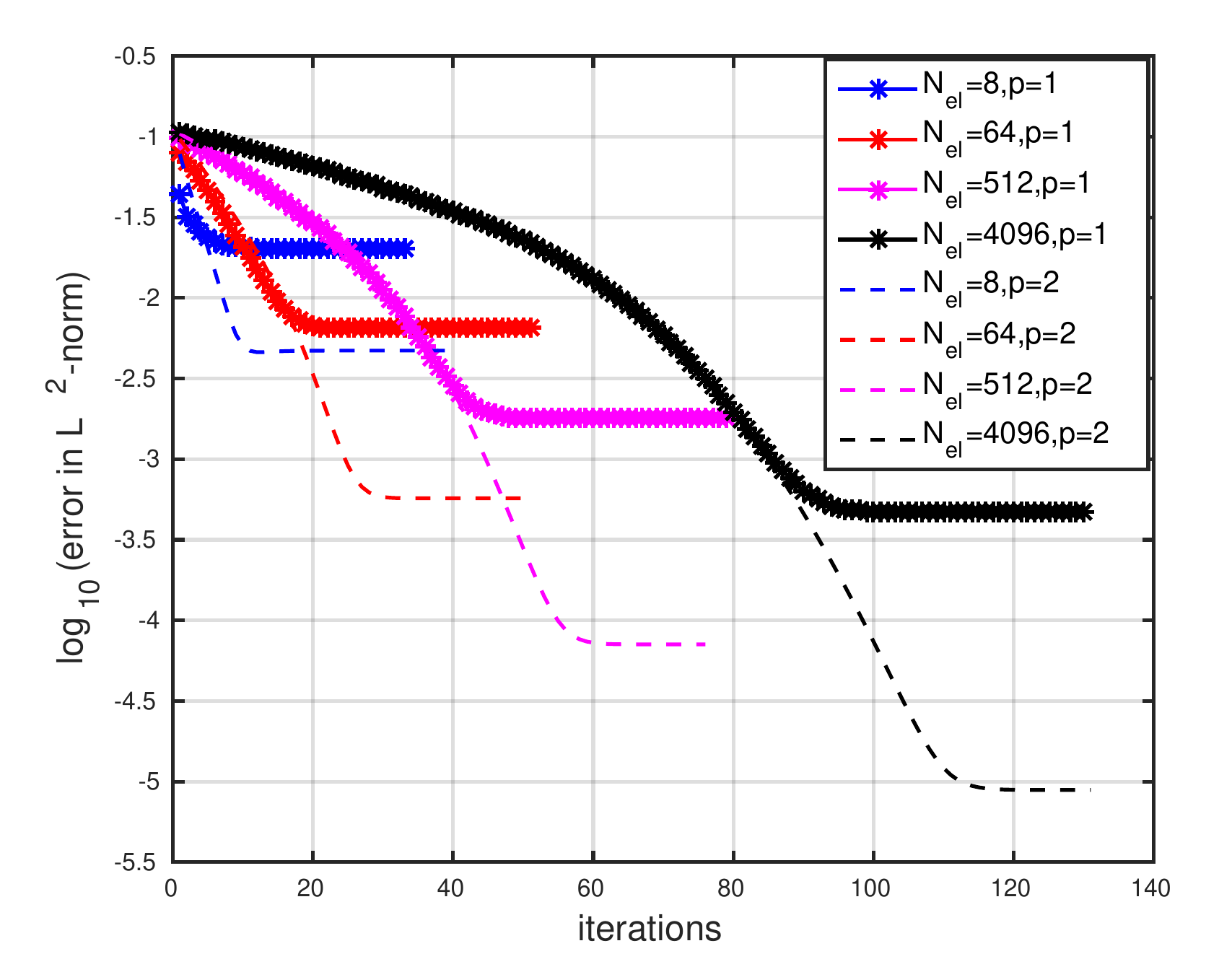}
}
\subfigure[Error history for p=3,4]{
\includegraphics[trim=2cm 7cm 2.5cm 7.25cm,clip=true,width=0.5\textwidth]{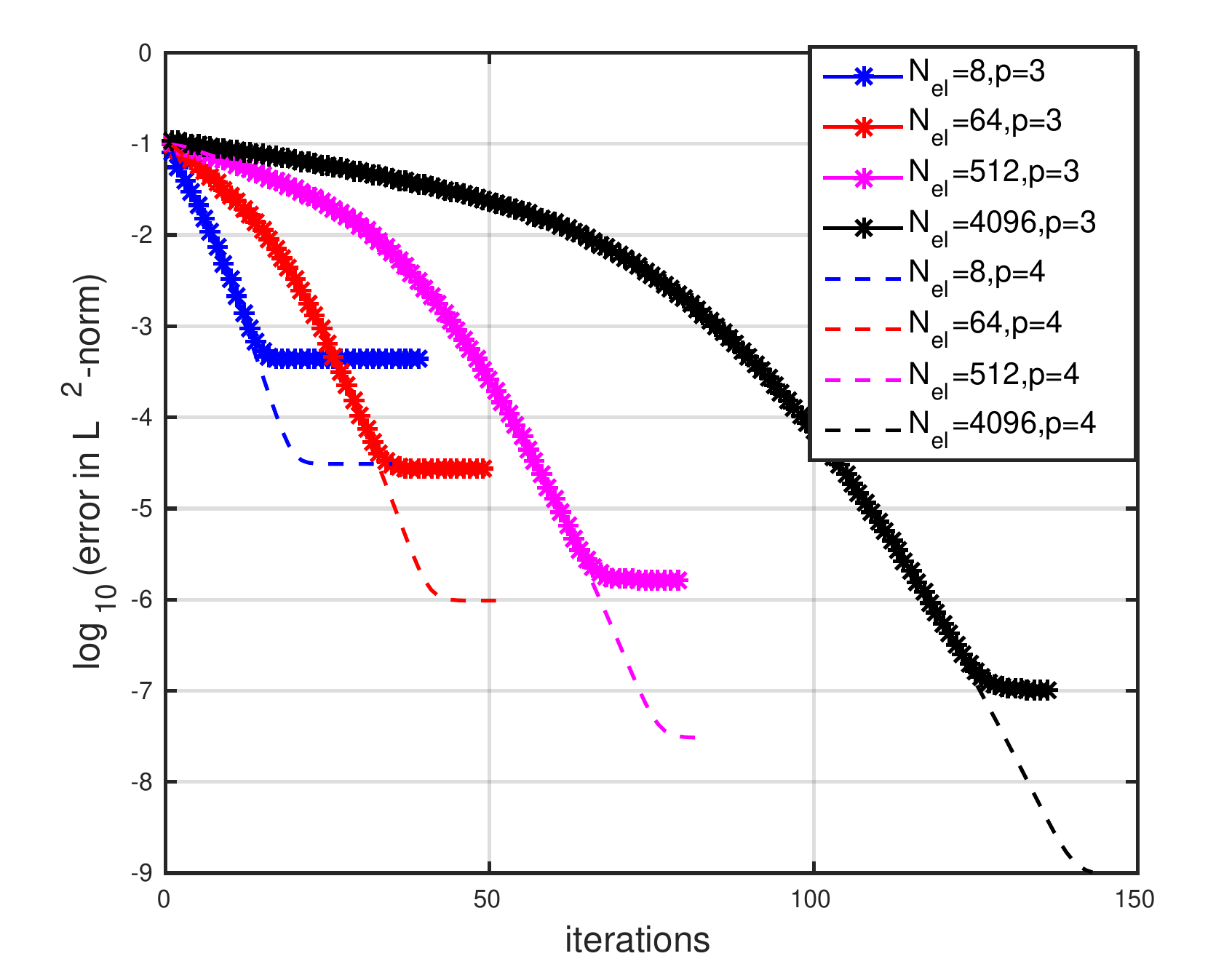}
}
\vspace{-5mm}
\caption{Convergence of eHDG for different $h$ and $p$ for 3D transport.} 
\label{allp_3dsteady}
\end{figure}
\vspace{-6mm}
\begin{figure}[h!t!b!]
  \subfigure[$\u$ at $iteration = 1$]{
    \includegraphics[trim=4cm 0cm 3cm 0cm,clip=true,width=0.48\columnwidth]{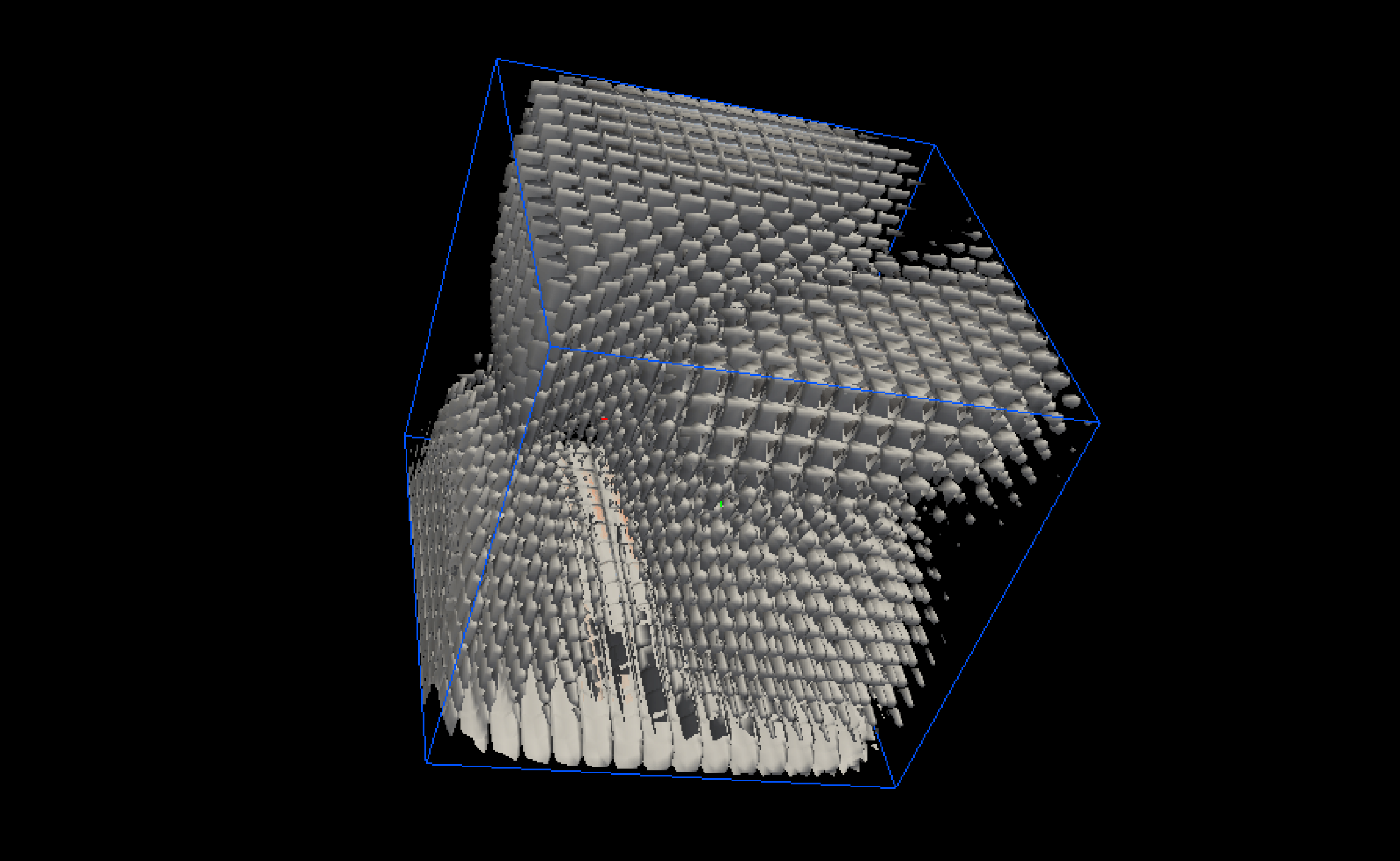}
  }
  \subfigure[$\u$ at $iteration = 16$]{
    \includegraphics[trim=4cm 0cm 3cm 0cm,clip=true,width=0.48\columnwidth]{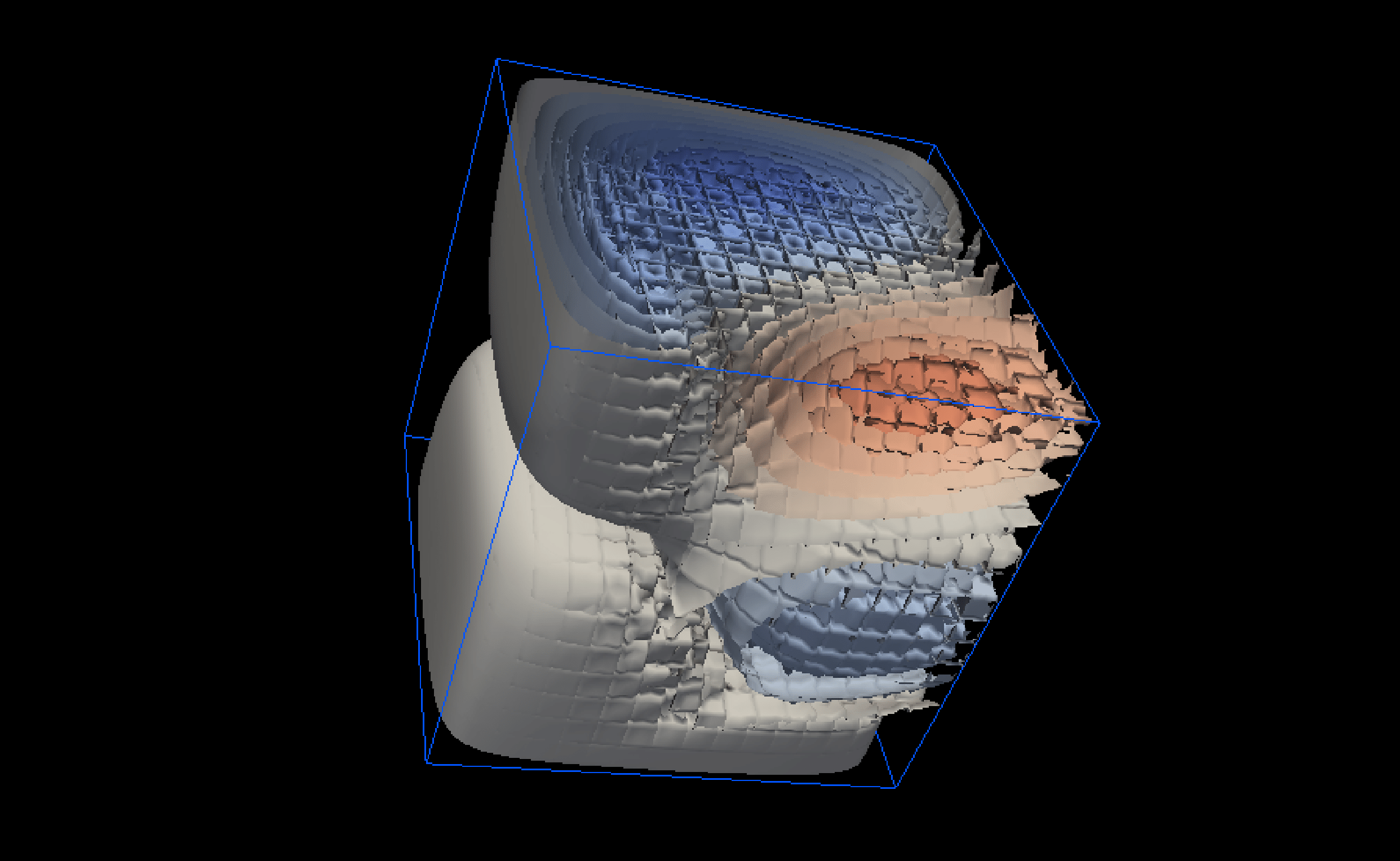}
  }
  \subfigure[$\u$ at $iteration = 48$]{
    \includegraphics[trim=4cm 0cm 3cm 0cm,clip=true,width=0.48\columnwidth]{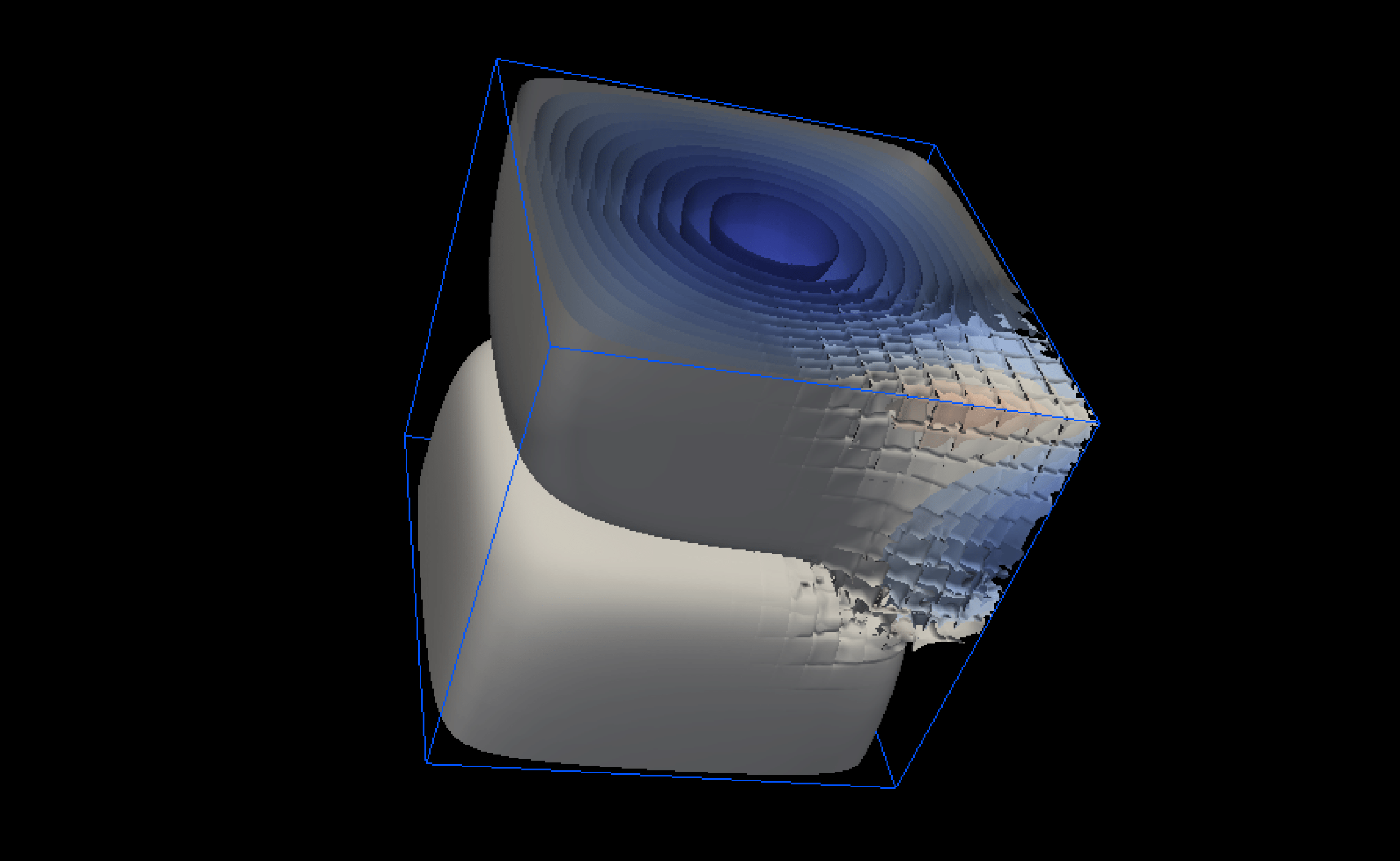}
  }
  \subfigure[$\u$ at $iteration = 143$]{
    \includegraphics[trim=4cm 0cm 3cm 0cm,clip=true,width=0.48\columnwidth]{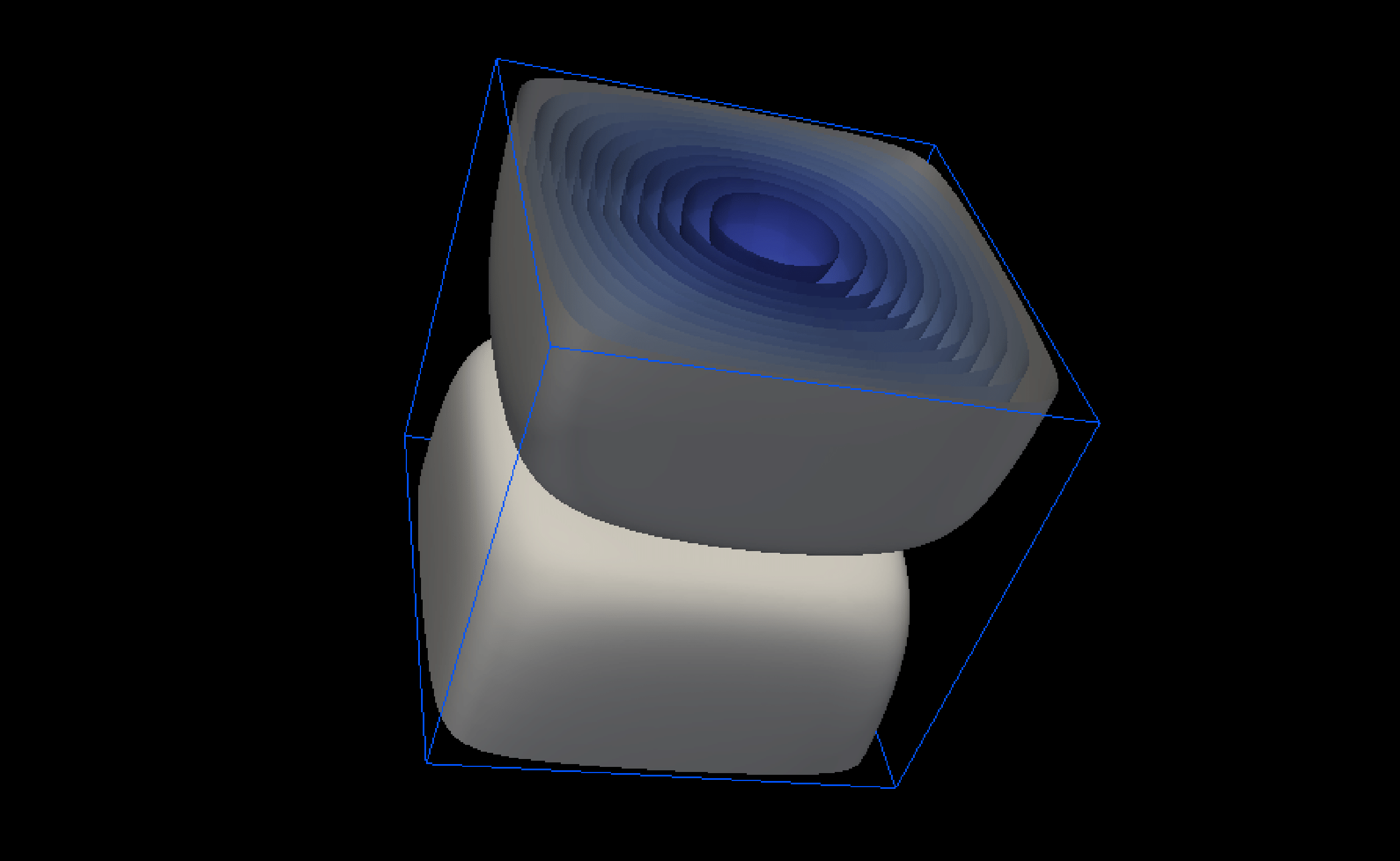}
  }
 \caption{Evolution of iterative eHDG solutions for 3D steady state transport equation.}
  \figlab{3diterfig1}
\end{figure}

The convergence history is shown in figure \ref{allp_3dsteady} and
they exhibit a similar trend as in 2D, i.e. exponential convergence
independent of the polynomial order (see also the $6$th column of
Table \ref{tab:Comparison_2dall_3d}. Again the iteration increases as
the mesh is refined.  The evolution of the eHDG solution with respect to
iterations in Figure \figref{3diterfig1} shows the convergence of
solution from inflow to outflow. Here, the solution order  is $p = 4$.

\subsection{2D linearized shallow water equations}
Here we consider equation \eqnref{linearizedShallow}, and in that we are considering a linear standing wave, which is an oceanic flow. 
For linear standing wave we take $\Phi=g=1$, $f=0$, $\gamma=0$ (zero bottom friction), $\taub=0$ (zero wind stress). The domain
is $[0,1]\times[0,1]$ and wall boundary condition is applied on the boundary. 
The following exact
solution \cite{GiraldoWarburton08} is taken
\begin{subequations}
\eqnlab{shallowexact}
\begin{align}
	\phi&=\cos(\pi x)\cos(\pi y)\cos(\sqrt{2} \pi t), \\
	u&=\frac{1}{\sqrt{2}}\sin(\pi x)\cos(\pi y)\sin(\sqrt{2} \pi t), \\
	v&=\frac{1}{\sqrt{2}}\cos(\pi x)\sin(\pi y)\sin(\sqrt{2} \pi t). 
\end{align}
\end{subequations}
The convergence of the $L^2$ norm of the solution is presented in
Figure \ref{l2err_comb}. Here we have taken $\Delta t=10^{-6}$ and
$10^5$ time steps in order to show the theoretical convergence rates
and from Figure \ref{l2err_comb-c} we see that optimal convergence
rate is obtained. The number of iterations required per time step in
this case is constant and is always equal to $2$ for all meshes and
polynomial orders considered. The reason is that the initial guess for
each time step is taken as the solution in the previous time step. Furthermore, 
the time step is small.

\begin{minipage}{\linewidth}
\centering
\captionof{table}{\it{Iterations per time step for 2D linearized shallow water equation and 3D time dependent advection for different $\Delta t$}} 
\label{tab:Comparison_shallow_water_3dadvec}
\begin{tabular}{ | r || r || r || c | c | c | c | }
\hline
\multirow{2}{*}{$N_{el}(2D)$} & \multirow{2}{*}{$N_{el}(3D)$} & \multirow{2}{*}{p} & \multicolumn{2}{|c|}{2D Shallow water} & \multicolumn{2}{|c|}{3D advection} \\
 & & & $\Delta t=10^{-3}$ & $\Delta t=10^{-4}$ & $\Delta t=10^{-3}$ & $\Delta t=10^{-4}$ \\
\hline
16 & 8 & 1 & 3 & 2 & 2 & 2 \\
\hline
64 & 64 & 1 & 4 & 2 & 3 & 2 \\
\hline
256 & 512 & 1 & 4 & 3 & 3 & 2 \\
\hline
1024 & 4096 & 1 & 4 & 3 & 3 & 2 \\
\hline
\hline
16 & 8 & 2 & 4 & 2 & 2 & 2 \\
\hline
64 & 64 & 2 & 4 & 2 & 3 & 2 \\
\hline
256 & 512 & 2 & 5 & 2 & 3 & 2 \\
\hline
1024 & 4096 & 2 & 6 & 3 & 3 & 2 \\
\hline
\hline
16 & 8 & 3 & 4 & 2 & 3 & 2 \\
\hline
64 & 64 & 3 & 4 & 2 & 3 & 2 \\
\hline
256 & 512 & 3 & 5 & 3 & 3 & 2 \\
\hline
1024 & 4096 & 3 & 6 & 3 & 4 & 2 \\
\hline
\hline
16 & 8 & 4 & 4 & 2 & 3 & 2 \\
\hline
64 & 64 & 4 & 5 & 3 & 3 & 2 \\
\hline
256 & 512 & 4 & 6 & 3 & 3 & 2 \\
\hline
1024 & 4096 & 4 & 7 & 3 & 4 & 2 \\
\hline
\end{tabular}
\end{minipage}
\vspace{5mm}

To compare with the 3D time-dependent advection in the next section,
we choose the time step of $\Delta t=10^{-3}$ and $\Delta t=10^{-4}$,
and tabulate the number of eHDG iterations in Table
\ref{tab:Comparison_shallow_water_3dadvec}. As can be seen, the number
of iterations increases slightly as we increase the solution order, and
this is consistent with Theorem \ref{ThmConvergenceShallowWater}.

\subsection{3D time dependent transport equation}
In this section we consider the following time-dependent transport equation
\begin{equation}
\eqnlab{timetransport}
\frac{\partial u}{\partial t}+\Div (\betab u)=f,
\end{equation}
and the exact solution is a Gaussian moving across the diagonal of a unit cube, i.e.,
\[
	u^e=e^{-5((x-0.2t)^2+(y-0.2t)^2+(z-0.2t)^2)},
\]
Structured hexahedral mesh $8\times8\times8$ is used and the solution order is $p = 4$. 
The time step is chosen $\Delta t=0.01$ and the simulation is run for
$240$ time steps.  Figure \figref{3ditertime} compares the numerical
solution using the eHDG iterative solver and the exact solution. The
tolerance criteria the same as in Section \ref{ex1_2dsteady}, and the
solver always takes 9 iterations per time step. In table
\ref{tab:Comparison_shallow_water_3dadvec} we compare the iterations
per time step required to converge for two smaller time step
sizes. Unlike shallow water equation, transport equation with eHDG
iterative solver has constant eHDG iterations as the solution
increases, and this is consistent with our theoretical result in Theorem \ref{ThmExponentialDecay}.

\begin{figure}[h!t!b!]
  \subfigure[$u$ at $time=0.5$]{
    \includegraphics[trim=4cm 0cm 3cm 0cm,clip=true,width=0.31\columnwidth]{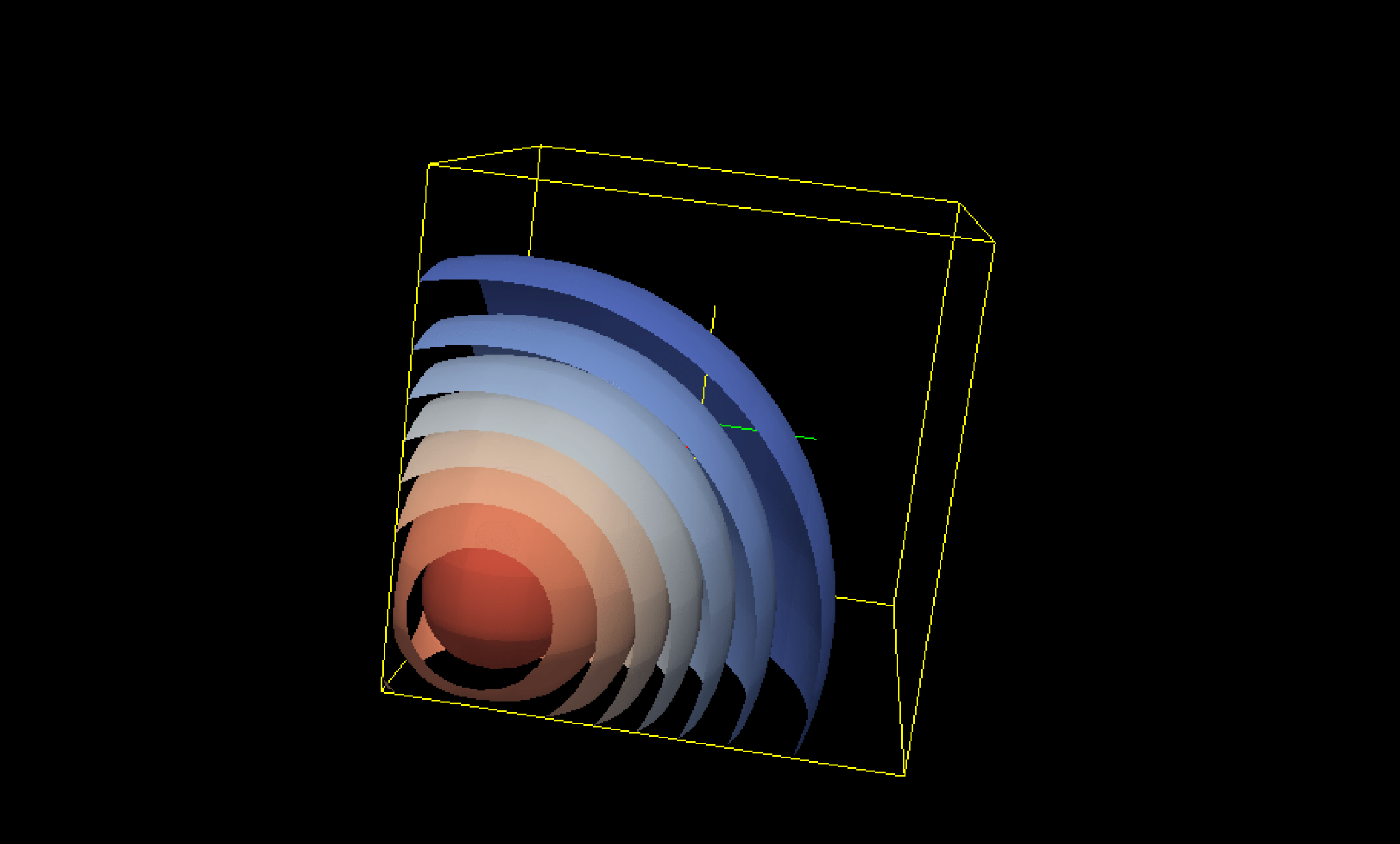}
  }
  \subfigure[$u$ at $time=1.4$]{
    \includegraphics[trim=4cm 0cm 3cm 0cm,clip=true,width=0.31\columnwidth]{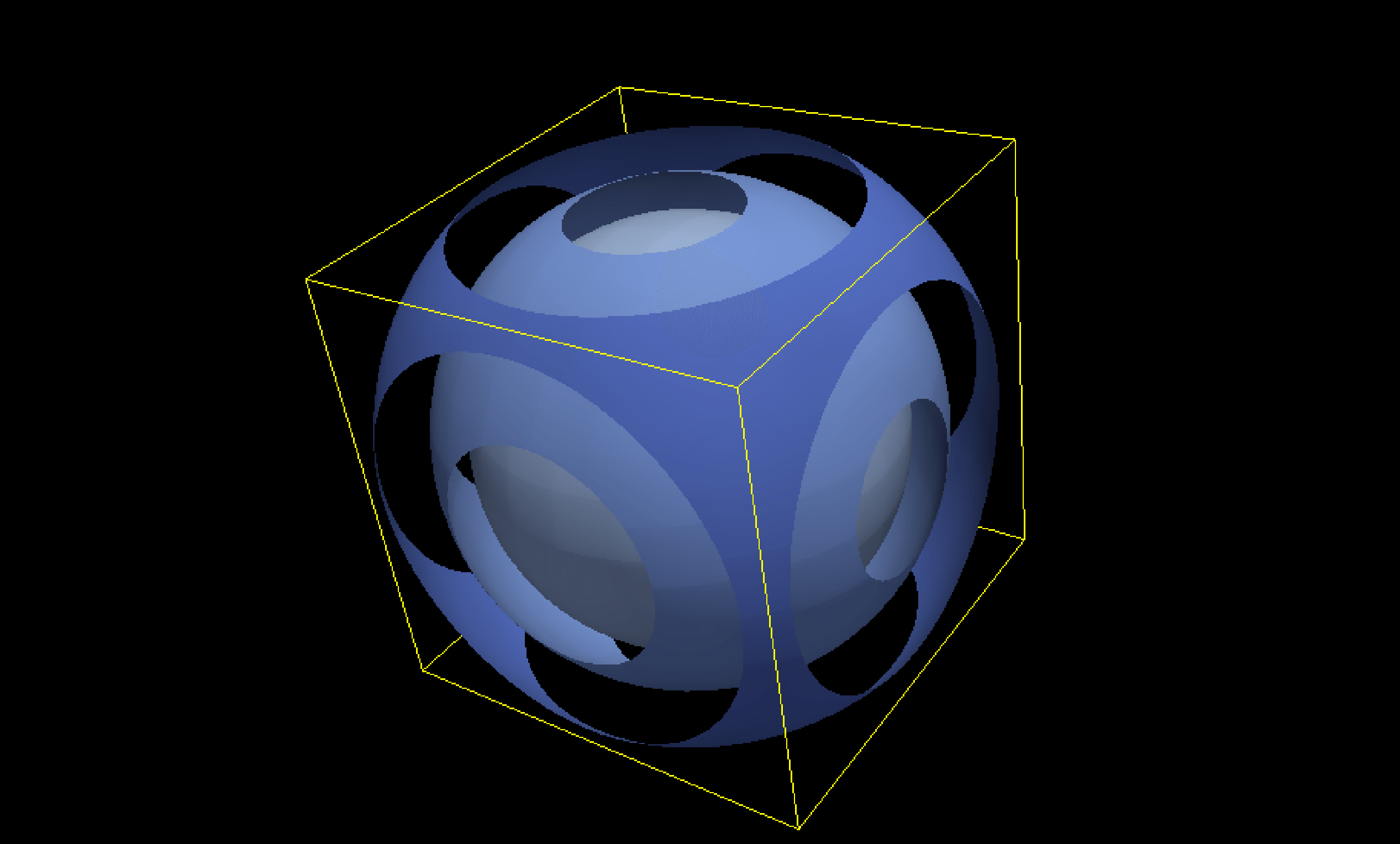}
  }
  \subfigure[$u$  at $time=2.4$]{
    \includegraphics[trim=4cm 0cm 3cm 0cm,clip=true,width=0.31\columnwidth]{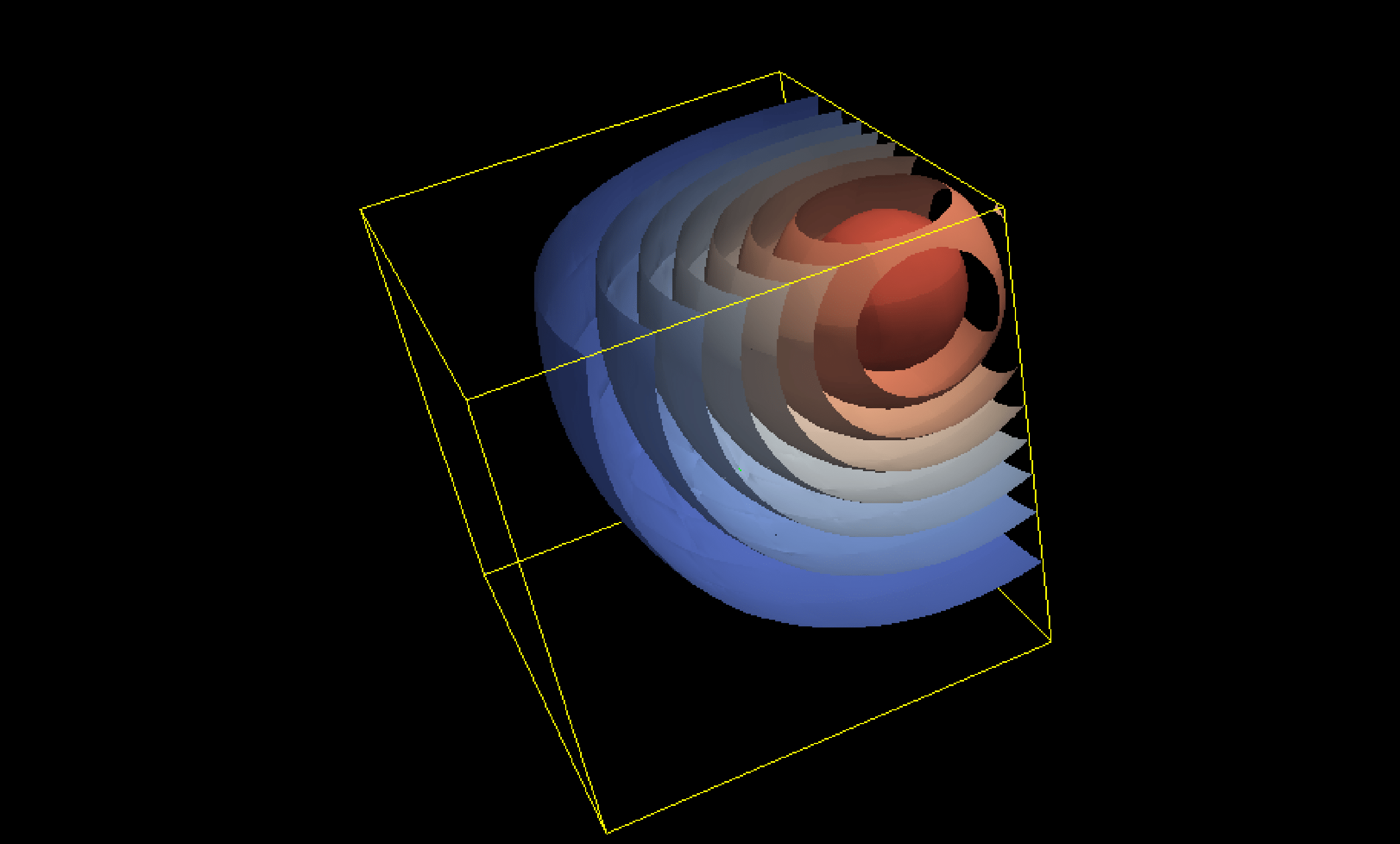}
  }
 \caption{eHDG solution for 3D time dependent transport equation}
  \figlab{3ditertime}
\end{figure}

\section{Conclusion}

We have presented an iterative solver, namely eHDG, for HDG
discretizations of hyperbolic systems. The method exploits the
structure of HDG discretization and idea from domain decomposition
methods.  The key features of the eHDG solver are: 1) it solves
independent element-by-element local equations during each iteration,
2) the number of iterations are independent of polynomial order, and 3) it
achieves exponential convergence rate. These features make the eHDG
solver naturally suitable for higher order HDG methods in large
scale parallel environments.



\bibliography{references,ceo}

\end{document}